\documentclass[12pt]{amsart}
\usepackage{amscd}
\usepackage{amssymb}

\usepackage{eepic}

\allowdisplaybreaks[1]

\newtheorem{thm}{Th\'eor\`eme}[section]

\newtheorem{prop}[thm]{Proposition}

\newtheorem{defn}[thm]{D\'efinition}
\theoremstyle{definition}

\numberwithin{equation}{section}
\newcommand{\R}{\mathbb R}

\newcommand{\n}{\noindent}
\renewcommand{\a}{\alpha}

\renewcommand{\b}{\beta}

\newcommand{\im}{{\rm im}}
\newcommand{\Hom}{{\rm Hom}}

\newcommand{\resumename}{R\'esum\'e}
\newenvironment{resume}{\narrower\footnotesize\bf
\noindent\resumename.\quad\footnotesize\rm}{\par\bigskip}

\font\hb=cmbx12
\newcommand\pt{\hbox{\hb .}}
\font\hbb=cmbx36
\newcommand\ppt{\hbox{\hbb .}}

\begin{document}

\title[Alg\`ebres \`a homotopie pr\`es, (co)homologies]{Alg\`ebres enveloppantes \`a homotopie pr\`es, homologies et cohomologies}
\author[R. Chatbouri]{Ridha Chatbouri}
\address{
D\'epartement de Math\'ematiques\\
Unit\'e de Recherche Physique Math\'ematique\\
Fa\-cult\'e des Sciences de Monastir\\
Avenue de l'environnement\\
5019 Monastir\\
Tuni\-sie} \email{Ridha.Chatbouri@ipeim.rnu.tn}

\maketitle


\begin{resume}
On pr\'esente une d\'efinition et une construction unif\'ee des homologies et cohomologies d'alg\`ebres et de modules sur ces alg\`ebres dans le cas d'alg\`ebres associatives ou commutatives ou de Lie ou de Gertsenhaber. On s\'epare la construction `lin\'eaire' des cog\`ebres ou bicog\`ebres qui traduisent les sym\'etries des relations de d\'efinition de la structure de la partie structure qui appara\^\i t ici comme une cod\'erivation de degr\'e 1 et de carr\'e nul de la cog\`ebre ou de la bicog\`ebre.\\

\end{resume}


\section{Introduction et motivation}\label{sec0}\

La d\'efinition op\'eradique de l'homologie permet d'unifier les notions d'homologie et de cohomologie pour les modules sur les alg\`ebres associatives (homologie et cohomologie de Hochschild) et pour les modules sur les alg\`ebres de Lie (homologie et cohomologie de Chevalley).\\

Ces deux complexes sont d\'eriv\'es de r\'esolutions de l'alg\`ebre (associative ou de Lie) consid\'er\'ee : la bar r\'esolution et la r\'esolution de Koszul.\\

D'autre part, on peut affaiblir la notion d'alg\`ebre pour d\'efinir des alg\`ebres `\`a homotopie pr\`es'. Dans le cas d'une alg\`ebre associative $A$, on regarde $A$ comme une alg\`ebre \`a homotopie pr\`es particuli\`ere : cette alg\`ebre \`a homotopie pr\`es est exactement la Bar r\'esolution de l'alg\`ebre $A$. Dans le cas d'une alg\`ebre de Lie $\mathfrak g$, la r\'esolution de Koszul est obtenue en tensorisant par $\mathcal U(\mathfrak g)$ l'alg\`ebre \`a homotopie pr\`es qu'on identifie \`a $\mathfrak g$.\\

Dans cet article, on parlera plut\^ot d'alg\`ebre enveloppante : alg\`ebre $A_\infty$ enveloppante de $A$, alg\`ebre $L_\infty$ enveloppante de $\mathfrak g$.\\

On d\'ecrira alors les liens entre d'une part l'homologie et la cohomologie d'un mo\-dule $M$ sur l'alg\`ebre  de type $X$ et des structures de $X_\infty$ alg\`ebre enveloppante associ\'ee.\\

De plus on traitera de m\^eme le cas des alg\`ebres commutatives $R$ et de leurs mo\-dules en \'etudiant leurs alg\`ebres \`a homotopie pr\`es enveloppantes $C_{\infty}(R)$. Enfin, nous utiliserons cette construction pour d\'efinir l'homologie et la cohomologie d'un module sur une alg\`ebre de Gerstenhaber $G$ gr\^ace \`a son alg\`ebre \`a homotopie pr\`es enveloppante.\\

Dans toutes ces constructions, nous commencerons par d\'efinir la partie lin\'eaire (les espaces de tenseurs retenus sur l'alg\`ebre de d\'epart) et sa ou ses comultiplications naturelles soit coassociatives $\Delta$ pour $A_\infty(A)$ ou $L_\infty(\mathfrak g)$, soit cocrochet de Lie $\delta$ pour $C_\infty(R)$ soit un couple ($\kappa,\Delta$) pour les structures d'alg\`ebres de Gerstenhaber \`a homotopie pr\`es $G_\infty(G)$. Ensuite nous \'ecrirons la structure de l'alg\`ebre $X_\infty$ comme l'op\'erateur bord d'un complexe. Cela nous permettra de d\'efinir explicitement les op\'erateurs d'homologie et de cohomologie associ\'es sur les modules.\\

Nous retrouvons ainsi des r\'esultats ant\'erieurs (voir \cite{[G]}, \cite{[BGHHW]} et \cite{[AAC3]} par exemple) en les compl\'etant, les pr\'ecisant et les explicitant.\\


\section{Homologie et cohomologie de Hochschild d'une alg\`ebre associative}\label{sec1}\

Soit $A$ une alg\`ebre associative et unitaire sur $\R$. On rappelle (voir par exemple \cite{[L]}) que la Bar-r\'esolution de $A$ est le complexe suivant~:
$$
0\stackrel{\partial_0}{\longleftarrow}A\stackrel{\partial_1}{\longleftarrow}A\otimes A\stackrel{\partial_2}{\longleftarrow}\dots \stackrel{\partial_{n-1}}{\longleftarrow}\otimes^nA\stackrel{\partial_n}{\longleftarrow}\dots
$$
o\`u $\partial_n$ est l'application lin\'eaire d\'efinie par~:
$$
\partial_n(a_0\otimes a_1\otimes\dots\otimes a_n)=\sum_{j=0}^{n-1}(-1)^ja_0\otimes\dots\otimes(a_ja_{j+1})\otimes\dots\otimes a_n.
$$

On v\'erifie directement que $\partial_{n-1}\circ\partial_n=0$ pour tout $n\geq1$. De plus, puisque $A$ est unitaire, ce complexe admet une homotopie, d\'efinie par~:
$$
h_n~:\otimes^nA\longrightarrow\otimes^{n+1}A,\qquad h_n(a_1\otimes a_2\otimes\dots\otimes a_n)=1\otimes a_1\otimes\dots\otimes a_n.
$$

On montre directement que $h_n\circ\partial_n+\partial_{n+1}\circ h_{n+1}$ est l'identit\'e de $\otimes^{n+1}A$. Ce complexe est donc une r\'esolution de l'alg\`ebre $A$. En posant $C_n(A)=\otimes^{n+1}A$, on le note~:
$$
0\stackrel{\partial_0}{\longleftarrow}C_0(A)\stackrel{\partial_1}{\longleftarrow}C_1(A)\stackrel{\partial_2}{\longleftarrow}\dots \stackrel{\partial_{n-1}}{\longleftarrow}C_{n-1}(A)\stackrel{\partial_n}{\longleftarrow}\dots
$$
On note $Z_n(A)=\ker(\partial_n)$ et $B_n(A)=\im(\partial_{n+1})$, on a donc
$$
H_n(A)=Z_n(A)/B_n(A)=0
$$
pour tout $n$.\\

Soit $M$ un $A$-bimodule. L'homologie de Hochschild du bimodule $M$ est d\'efinie \`a partir de la Bar-r\'esolution de la fa\c con suivante.

On pose $C_n(A,M)=\sum_{j=0}^n\otimes^jA\otimes M\otimes\otimes^{n-j}A$. Le complexe d'homologie de Hochschild est le complexe~:
$$
0\stackrel{\partial_0}{\longleftarrow}C_0(A,M)\stackrel{\partial_1}{\longleftarrow}C_1(A,M)\stackrel{\partial_2}{\longleftarrow}\dots \stackrel{\partial_{n-1}}{\longleftarrow}C_{n-1}(A,M)\stackrel{\partial_n}{\longleftarrow}\dots
$$
o\`u le bord $\partial_n$ est lin\'eaire et d\'efini par~:
$$\aligned
\partial_n(a_0\otimes\dots&\otimes a_{i-1}\otimes m_i\otimes a_{i+1}\otimes\dots\otimes a_n)=\cr
&\sum_{j=0}^{i-2}(-1)^ja_0\otimes\dots\otimes a_ja_{j+1}\otimes\dots\otimes m_i\otimes\dots\otimes a_n+\cr
&+(-1)^{i-1}a_0\otimes\dots\otimes a_{i-1}m_i\otimes\dots\otimes a_n+(-1)^ia_0\otimes\dots\otimes m_ia_{i+1}\otimes\dots\otimes a_n+\cr
&+\sum_{j=i+1}^{n-1}(-1)^ja_0\otimes\dots\otimes m_i\otimes\dots\otimes a_ja_{j+1}\otimes\dots\otimes a_n
\endaligned
$$
Comme plus haut, le $n^{i\grave{e}me}$ groupe d'homologie du bimodule $M$ est $H_n(A,M)=Z_n(A,M)/B_n(A,M)$, o\`u $Z_n(A,M)=\ker(\partial_n)$ et $B_n(A,M)=\im(\partial_{n+1})$.\\

De m\^eme la cohomologie de ce bimodule $M$ est d\'efinie ainsi. On pose
$$
C^n(A,M)=L(\otimes^nA,M),
$$
et on obtient le complexe de cohomologie de Hochschild d\'eriv\'e de la Bar r\'esolution~:
$$
0\stackrel{\partial^0}{\longrightarrow}C^0(A,M)\stackrel{\partial^1}{\longrightarrow}C^1(A,M)\stackrel{\partial^2}{\longrightarrow}\dots \stackrel{\partial^{n-1}}{\longrightarrow}C^{n-1}(A,M)\stackrel{\partial^n}{\longrightarrow}\dots
$$
avec~:
$$\aligned
(\partial^{n+1}f_n)(a_0\otimes\dots\otimes a_n)&=a_0f_n(a_1\otimes\dots\otimes a_n)-\sum_{j=0}^{n-1}(-1)^jf_n(a_0\otimes\dots\otimes a_ja_{j+1}\otimes\dots\otimes a_n)\cr
&\hskip 2cm+(-1)^{n-1}f_n(a_0\otimes\dots\otimes a_{n-1})a_n.
\endaligned
$$
En consid\'erant l'isomorphisme $L(\otimes^nA,M)\simeq\Hom_{A\otimes A^{op}}(\otimes^{n+2}A,M)$, donn\'e par
$$
F_n(a_0\otimes\dots\otimes a_{n+1})=a_0 f_n(a_1\otimes\dots\otimes a_n)a_{n+1},
$$
$(A^{op}$ est $A$ munie du produit oppos\'e $(a,b)\longmapsto ba$), on constate qu'on peut \'ecrire
$$
\partial^{n+1} f_n=-F_n\circ\partial_{n+1}|_{1\otimes(\otimes^{n+1}A)\otimes1},
$$
donc on a un complexe d\'eriv\'e de la Bar r\'esolution. La cohomologie de Hochschild de $M$ est celle de ce complexe~:
$$
H^n(A,M)=Z^n(A,M)/B^n(A,M),\text{ o\`u } Z^n(A,M)=\ker(\partial^{n+1})\text{ et } B^n(A,M)=\im(\partial^n).
$$

\section{Alg\`ebre $A_\infty$ enveloppante}\label{sec2}\

Dans cette section, nous allons traduire les constructions pr\'ec\'edentes dans le cadre de la $A_\infty$ alg\`ebre enveloppante de $A$.\\

La structure multiplicative de $A$ n'\'etant pas commutative, ne poss\`ede pas de sym\'etrie~: l'unique relation est l'associativit\'e~:
$$
Ass(a_1,a_2,a_3)=a_1(a_2a_3)-(a_1a_2)a_3,
$$
L'application $Ass$ n'a pas de sym\'etries, l'espace $Vect(Ass\circ\sigma,~~\sigma\in\mathfrak{S}_3)$ est de dimension 6. On consid\'erera donc toute l'alg\`ebre tensorielle de $A$.\\

D'autre part, on sait que les signes $(-1)^j$ apparaissant dans les formules ci-dessus s'interpr\`etent simplement si on d\'ecale le degr\'e des \'el\'ements de $A$ de 1. On va donc directement consid\'erer ici une alg\`ebre $A$ gradu\'ee, munie d'une multiplication de degr\'e 0~:
$$
|a_1a_2|=|a_1|+|a_2|,
$$
si $|a|$ est le degr\'e de $a\in A$.\\

L'espace vectoriel que l'on consid\`ere est donc $\oplus_{n\geq1}\otimes^nA[1]$, l'alg\`ebre tensorielle sans unit\'e, gradu\'ee par~:
$$
\deg(a_1\otimes\dots\otimes a_n)=\deg(a_1)+\dots+\deg(a_n)=|a_1|+\dots+|a_n|-n.
$$
Cet espace est une cog\`ebre coassociative libre pour la comultiplication canonique~:
$$
\Delta(a)=0,\quad\Delta(a_1\otimes\dots\otimes a_n)=\sum_{j=1}^{n-1}(a_1\otimes\dots\otimes a_j)\bigotimes(a_{j+1}\otimes\dots\otimes a_n).
$$
(Voir par exemple \cite{[AAC3]}). La coassociativit\'e s'\'ecrit~:
$$
(id\otimes\Delta)\circ\Delta=(\Delta\otimes id)\circ\Delta.
$$
Le fait que $(\oplus_{n\geq1}\otimes^nA[1],\Delta)$ est une cog\`ebre libre se traduit ici par la propri\'et\'e suivante~:
\begin{itemize}
\item[1.] Pour tout espace gradu\'e $V$, un morphisme $F$ de cog\`ebres de $(\oplus_{n\geq1}\otimes^nA[1],\Delta)$ dans $(\oplus_{n\geq1}\otimes^nV[1],\Delta)$ est uniquement caract\'eris\'e par sa projection sur $V[1]$~: $F$ est un morphisme si $(F\otimes F)\circ \Delta=\Delta\circ F$. Posons $F_n~:\otimes^nA[1]\longrightarrow V[1]$, on a explicitement~:
$$
F(a_1\otimes\dots\otimes a_n)=\sum_{k=1}^n\hskip 0.1cm\sum_{r_j,~r_1+\dots+r_k=n}\hskip 0.1cm F_{r_1}(a_1\otimes\dots\otimes a_{r_1})\otimes \dots\otimes F_{r_k}(a_{n-r_k+1}\otimes\dots\otimes a_n).
$$
R\'eciproquement, pour toute suite d'applications lin\'eaires $(F_n)$, l'application $F$ d\'efinie ci-dessus est un morphisme de cog\`ebres.\\

\item[2.] Toute cod\'erivation $D$ de $(\oplus_{n\geq1}\otimes^nA[1],\Delta)$ est uniquement caract\'eris\'ee par sa projection sur $A[1]$~: $D$ est une d\'erivation si $(D\otimes id+id\otimes D)\circ\Delta=\Delta\circ D$. Posons $D_n~:\otimes^nA[1]\longrightarrow A[1]$, on a explicitement~:
$$\aligned
D(a_1\otimes\dots\otimes a_n)&=\sum_{1\leq j<j+r-1\leq n}\hskip 0.1cm(-1)^{\deg(D_r)(\deg(a_1)+\cdots+deg(a_{j-1}))}\cr
&(a_1\otimes\dots\otimes a_{j-1})\otimes D_r(a_j\otimes\dots\otimes a_{j+r-1})\otimes(a_{j+r}\otimes\dots\otimes a_n).
\endaligned
$$
R\'eciproquement, pour toute suite d'applications lin\'eaires $(D_n)$, l'application $D$ d\'efinie ci-dessus est une cod\'erivation de la cog\`ebre $(\oplus_{n\geq1}\otimes^nA[1],\Delta)$.\\

\end{itemize}

La multiplication devient une op\'eration de degr\'e 1 de $\otimes^2A[1]$ dans $A[1]$. On utilise la r\`egle des signes donn\'ee dans \cite{[AAC2]}, et on lui associe l'application~:
$$
m(a_1\otimes a_2)=(-1)^{\deg(a_1)}(a_1a_2).
$$

Alors, l'application $m$ se prolonge d'une fa\c con unique en une cod\'erivation de degr\'e 1 de la loi $\Delta$ \`a tout l'espace $\oplus_{n\geq1}\otimes^nA[1]$.

Ce prolongement est explicitement~:
$$
m(a_0\otimes\dots\otimes a_n)=\sum_{j=0}^{n-1}(-1)^{\sum_{i<j}\deg(a_i)}a_0\otimes\dots\otimes m(a_j\otimes a_{j+1})\otimes\dots\otimes a_n.
$$

La relation d'associativit\'e est \'equivalente \`a $m\circ m=0$ (\'equation de structure).\\

\begin{defn}{\rm(Alg\`ebre $A_\infty$ enveloppante)}

\

Soit $A$ une alg\`ebre associative. On appelle alg\`ebre $A_\infty$ enveloppante de $A$ la cog\`ebre $(\oplus_{n\geq1}\otimes^nA[1],\Delta)$ munie de la cod\'erivation $m$. On note cette alg\`ebre enveloppante
$$
A_\infty(A)=(\oplus_{n\geq1}\otimes^nA[1],\Delta,m).
$$
\end{defn}

La $A_\infty$ alg\`ebre enveloppante de $A$ est donc la g\'en\'eralisation au cas gradu\'e de la Bar r\'esolution de $A$~: si $A$ n'est pas gradu\'ee, tous les $a$ de $A[1]$ sont de degr\'e -1, et l'espace vectoriel $A_\infty(A)$ est $\oplus_{n\geq1}\otimes^nA$ et la cod\'erivation $m$ co\"\i ncide avec $\partial$.\\

Plus g\'en\'eralement,\\

\begin{prop}{\rm(Homologie et cohomologie de Hochschild et $A_\infty$ alg\`ebres)}

\

Soit $A$ une alg\`ebre associative, unitaire et gradu\'ee. Soit $M$ un $A$-bimodule gradu\'e. On d\'efinit une nouvelle alg\`ebre associative, not\'ee $B=A\ltimes M$, produit semi-direct de $A$ par $M$ en munissant $A\oplus M$ de la multiplication~:
$$
(a+u)(b+v)=(ab+(ub+av))\qquad (a,~b\in A,~~u,~v\in M).
$$
Alors
\begin{itemize}
\item[1.] Le complexe d'homologie de Hochschild du bimodule $M$ est un sous complexe de $(A_\infty(A\ltimes M),m_{A\ltimes M})$.\\

\item[2.] Un morphisme de cog\`ebres $F:(\oplus_{n\geq1}\otimes^nA[1],\Delta_A)\longrightarrow(\oplus_{n\geq1}\otimes^nB[1],\Delta_B)$, de degr\'e $0$ d\'efini par $F_1=\iota+c_1$, $F_n=c_n$ o\`u $\iota$ est l'injection canonique de $A$ dans $A\ltimes M$ et $c_j\in L(\otimes^jA[1],M[1])$ est un morphisme
$$
F~:(A_\infty(A),\Delta_A,m_A)\longrightarrow(A_\infty(B),\Delta_B,m_B)
$$
de $A_\infty$ alg\`ebres si et seulement si $m_B\circ F=F\circ m_A$, si et seulement si $\partial c_j=0$ ($j\geq1$) o\`u $\partial$ est le cobord de Hochschild.

\item[3.] Un tel morphisme est dit trivial s'il existe une suite d'applications $b=(b_n)$, de degr\'e -1 telle que $c=m_B\circ(\iota\otimes b-b\otimes\iota)-b\circ m_A$, ceci est \'equivalent à $c_1=0$ et $c_n=\partial b_{n-1}$, pour tout $n>1$.\\
\end{itemize}
\end{prop}

\

\noindent
{\bf Preuve}

\

\noindent
1. On a $B=A\oplus M$, alors, 
$$
C_n(B)=\otimes^{n-1}B\supset\sum_{j=0}^{n-1}\otimes^jA\otimes M\otimes^{n-1-j}A=C_n(A,M).
$$
On v\'erifie directement que ${m_B}_{|C_n(A,M)}=\partial_n$.

D'o\`u le complexe d'homologie de Hochschild du bimodule $M$ est
un sous complexe de $\Big(A_\infty(A\ltimes M), m_{A\ltimes M}\Big)$.

\noindent
2. Soit $F:\Big(\displaystyle\bigoplus_{n\geq1}\otimes^n
A[1],\Delta_A,m_A\Big)\longrightarrow\Big(\displaystyle\bigoplus_{n\geq1}\otimes^n
B[1],\Delta_B,m_B\Big)$ un morphisme de $A_\infty$-alg\`ebres tel
que $F_1=\iota+c_1$ et $F_n=c_n$. Alors, on a

$$F(a_1\otimes\dots \otimes a_n)=\sum_{\begin{smallmatrix}k>0,~0<r_1,\dots,r_k\\ r_1+\dots+r_k=n\end{smallmatrix}}F_{r_1}(a_1\otimes\dots \otimes a_{r_1})
\otimes \dots\otimes F_{r_k}(a_{n-r_k+1}\otimes\dots \otimes a_n).
$$ et $m_B\circ F=F\circ m_A$. On va d\'emontrer que $\partial
c_j=0$ par r\'ecurrence sur $j$.

Pour $j=1$, $F(a\otimes b)=F_1(a)\otimes F_1(b)+ F_2(a\otimes b)$.
Donc,
$$
m_B\circ F(a\otimes b)= m_B(F_1(a)\otimes
F_1(b)=(-1)^{deg a}(ab+c_1(a)b+ac_1(b))
$$
$$
=F\circ m_A(a\otimes b)= (-1)^{dega}(ab+c_1(ab))
.$$

D'o\`u $\partial c_1(a\otimes b)=0$.

Supposons que $\partial c_j=0$, $\forall j\leq n-2$ et montrons
que $\partial c_{n-1}=0$. On a
\begin{align*}
F\circ m_A&(a_0\otimes\dots\otimes a_n)=\cr
&=F(\sum_{j=0}^{n-1}(-1)^{\sum_{i<j}\deg(a_i)}a_0\otimes\dots\otimes m_A(a_j\otimes a_{j+1})\otimes\dots\otimes a_n)\\&=\sum_{\begin{smallmatrix}k>0,~0<r_1,\dots,r_k\\ r_1+\dots+r_k=n\end{smallmatrix}}\sum_{j=0}^{n-1}(-1)^{\sum_{i<r_1+...+r_{ j-1}+1}\deg(a_i)}F_{r_1}(a_1\otimes\dots \otimes a_{r_1})
\otimes\\
&\hskip 1cm\otimes\dots\otimes F_{r_j}\circ m_A(a_{r_1+...+r_{ j-1}+1}\otimes\dots \otimes a_{r_j})
\otimes \dots\otimes F_{r_k}(a_{n-r_k+1}\otimes\dots \otimes a_n)
\end{align*}
et
\begin{align*}
m_B \circ F&(a_0\otimes\dots\otimes a_n)=\\
=&\sum_{\begin{smallmatrix}k>0,~0<r_1,\dots,r_k\\ r_1+\dots+r_k=n\end{smallmatrix}}\sum_{j=0}^{n-1}(-1)^{\sum_{i<r_1+...+r_{ j-1}+1}\deg(a_i)}F_{r_1}(a_1\otimes\dots \otimes a_{r_1})
\otimes\dots\otimes \\&\otimes  m_B\Big( F_{r_j}(a_{r_1+...+r_{ j-1}+1}\otimes\dots \otimes a_{r_1+...+r_j})\otimes F_{r_{j+1}}(a_{r_1+...+r_{ j}+1}\otimes\dots \otimes a_{r_1+...+r_{j+1}})\Big)\\&
\otimes \dots\otimes F_{r_k}(a_{n-r_k+1}\otimes\dots \otimes a_n).\end{align*}

Dans la somme pr\'ec\'edente, si $r_j$ et $r_{j+1}$ sont sup\'erieurs \`a 1, on a
$$
m_B\Big( F_{r_j}(a_{r_1+...+r_{ j-1}+1}\otimes\dots \otimes a_{r_1+...+r_j})\otimes F_{r_{j+1}}(a_{r_1+...+r_{ j}+1}\otimes\dots \otimes a_{r_1+...+r_{j+1}})\Big)=0.
$$

De plus, on a si $r_j=1$,
$$\aligned
m_B&\Big(F_1(a_{r_1+...+r_{j-1}+1})\otimes F_{r_{j+1}}(a_{r_1+...+r_{j}+1}\otimes\dots \otimes a_{r_1+...+r_{j+1}})\Big)=\\
&\hskip 4cm=m_B\Big(a_{r_1+...+r_j}\otimes F_{r_{j+1}}(a_{r_1+...+r_j+1}\otimes\dots\otimes a_{r_1+...+r_{j+1}})\Big)
\endaligned
$$
et si $r_{j+1}=1$,
$$\aligned
m_B&\Big(F_{r_j}(a_{r_1+...+r_{j-1}+1}\otimes\dots \otimes a_{r_1+...+r_j})\otimes F_1(a_{r_1+...+r_j+1})\Big)=\\
&\hskip 4cm=m_B\Big(F_{r_j}(a_{r_1+...+r_{ j-1}+1}\otimes\dots \otimes a_{r_1+...+r_j})\otimes a_{r_1+...+r_j+1}\Big).
\endaligned
$$

Dans l'expression
$ (m_B \circ F-F\circ m_A)(a_0\otimes\dots\otimes a_n)$, les termes o\`u $r_j\leq n-2$ et les termes o\`u $r_{j+1}\leq n-2$  disparaissent gr\^ace \`a notre  hypoth\`ese de r\'ecurrence car $\partial c_{r_j}=0$, pour tout $ r_j\leq n-2$. Il ne reste que les termes o\`u $r_j=1$ et $r_{j+1}=n-1$ ou $r_j=n-1$ et $r_{j+1}=1$. Donc
$$
\begin{aligned}
(m_B \circ F&-F\circ m_A)(a_0\otimes\dots\otimes a_n)=m_B\big(F_{n-1}(a_0\otimes\dots\otimes a_{n-1})\otimes a_n\big)\\
&+m_B\big(a_0\otimes F_{n-1}(a_1\otimes\dots\otimes a_{n})\big)\cr
&-\sum_{j=1}^{n-1}(-1)^{\sum_{i<j}deg a_i}F_{n-1}\big(a_0\otimes\dots\otimes a_{j-1}\otimes m_A
(a_j\otimes a_{j+1})\otimes a_{j+2}\otimes\dots\otimes a_n\big)\cr
=~~&(-1)^{\sum_{i\leq n-1}deg a_i}F_{n-1}(a_0\otimes\dots\otimes a_{n-1}).a_n+(-1)^{deg a_0}a_0.F_{n-1}(a_1\otimes\dots\otimes a_{n})\cr
&-\sum_{j=1}^{n-1}(-1)^{\sum_{i\leq j}deg a_i}F_{n-1}\big(a_0\otimes\dots\otimes a_{j-1}\otimes
(a_j.a_{j+1})\otimes a_{j+2}\otimes\dots \otimes a_n\big)\cr
=~~&(\partial F_{n-1})(a_0\otimes\dots\otimes a_n)=0
\end{aligned}
$$
On retrouve l'op\'erateur de cobord de Hochschild $\partial $.

\

\noindent
3. Si $C=m_B\circ(\iota\otimes b-b\otimes\iota)-b\circ m_A$, alors,
$$
\begin{aligned}C_{n-1}(a_0\otimes&\dots\otimes a_n)=\\
=&m_B(a_0\otimes b_{n-2}(a_1\otimes\dots\otimes a_n)+b_{n-2}(a_0\otimes\dots a_{n-1})\otimes a_n)\\
&-\sum_{j=0}^{n-1}(-1)^{\sum_{i\leq j}deg a_i}b_{n-2}(a_0\otimes\dots \otimes a_ja_{j+1}\otimes \dots\otimes a_n)\\
=&(-1)^{deg a_0}a_0. b_{n-2}(a_1\otimes\dots\otimes a_n)-(-1)^{\sum_{i\leq n-1}deg a_i}b_{n-2}(a_0\otimes\dots \otimes a_{n-1}). a_n\\&-\sum_{j=0}^{n-1}(-1)^{\sum_{i\leq j}deg a_i}b_{n-2}(a_0\otimes\dots \otimes a_ja_{j+1}\otimes \dots\otimes a_n)\\
=&\partial b_{n-2}(a_0\otimes\dots\otimes a_n).
\end{aligned}$$

\

\vskip0.4cm

Dans les sections suivantes, nous allons r\'ep\'eter cette construction dans le cas des alg\`ebres commutatives, des alg\`ebres de Lie et des alg\`ebres de Gerstenhaber.\\


\section{Alg\`ebre $C_\infty$ enveloppante d'une alg\`ebre commutative}\label{sec3}\

Dans ce paragraphe, $R$ est une alg\`ebre associative, commutative, gradu\'ee sur $\R$. La relation de commutativit\'e~:
$$
Com(a_1\otimes a_2)=a_1a_2-(-1)^{|a_1||a_2|}a_2a_1=0
$$ 
est antisym\'etrique. La relation d'associativit\'e n'est pas compl\`etement antisym\'etrique mais on a~:
$$
\dim Vect\{Ass\circ\sigma,~~\sigma\in\mathfrak{S}_3\}=2.
$$
En fait$(Ass\circ\sigma)(a_1,a_2,a_3)=0$ pour tout $\sigma$ est \'equivalent, en tenant compte de la commutativit\'e, aux 2 relations~:
$$
Ass(a_1,a_2,a_3)=0\quad\text{ et }\quad Ass(a_2,a_1,a_3)=0.
$$
Un suppl\'ementaire de l'espace $Vect(id,(1,2))$ est engend\'e par les (2,1) et (1,2) battements de trois lettres.

On d\'efinit donc d'abord les $p,q$ battements de $n=p+q$ lettres comme les permutations $\sigma$ de $\{1,\dots,n\}$ telles que $\sigma(1)<\dots<\sigma(p)$ et $\sigma(p+1)<\dots<\sigma(p+q)$. On appelle $Bat(p,q)$ l'ensemble de tous ces battements et on d\'efinit le produit battement de deux tenseurs $\alpha=a_1\otimes\dots\otimes a_p\in\otimes^p R$ et $\beta=a_{p+1}\otimes\dots\otimes a_{p+q}\in\otimes^q R$ par
$$
bat_{p,q}(\alpha,\beta)=\sum_{\sigma\in Bat(p,q)}\varepsilon_{|a|}(\sigma^{-1}) a_{\sigma^{-1}(1)}\otimes\dots\otimes a_{\sigma^{-1}(n)}.
$$
O\`u $\varepsilon_{|a|}(\sigma^{-1})$ est le signe obtenu par la r\`egle de Koszul en tenant compte des degr\'es $|a|$ des tenseurs. C'est \`a dire le morphisme de groupe sur ${\mathfrak S}_n$ tel que $\varepsilon_{|a|}((i,j))=(-1)^{|a_i||a_j|}$.

$bat_{p,q}$ repr\'esente la somme sign\'ee de tous les tenseurs $a_{i_1}\otimes\dots\otimes a_{i_n}$ dans lesquels les vecteurs $a_1,\dots,a_p$ et les vecteurs $a_{p+1},\dots,a_{p+q}$ apparaissent rang\'es dans leur ordre naturel.\\

Par d\'efinition, l'espace $\underline{\otimes^nR}$ est le quotient de $\otimes^nR$ par la somme de toutes les images des applications lin\'eaires $bat_{p,n-p}$ ($0<p<n$) (voir \cite{[G]}, \cite{[L]}). Cet espace est engendr\'e par les classes des tenseurs $a_1\otimes\dots\otimes a_n$ dans le quotient. On note cette classe~: 
$$
\underline{a_1\otimes\dots\otimes a_n}.
$$

Si $n=3$, si $(a_1,a_2,a_3)$ est un syst\`eme libre dans $R$, l'espace $Bat(2,1)$ contient les tenseurs
$$\aligned
bat_{2,1}(a_1\otimes a_2,a_3)&=a_1\otimes a_2\otimes a_3+(-1)^{|a_3||a_2|}a_1\otimes a_3\otimes a_2+(-1)^{|a_3|(|a_2|+|a_1|)}a_3\otimes a_1\otimes a_2\cr
bat_{2,1}(a_2\otimes a_1,a_3)&=a_2\otimes a_1\otimes a_3+(-1)^{|a_3||a_1|}a_2\otimes a_3\otimes a_1+(-1)^{|a_3|(|a_2|+|a_1|)}a_3\otimes a_2\otimes a_1.
\endaligned
$$
L'espace $Bat(1,2)$ contient les tenseurs
$$\aligned
bat_{1,2}(a_1,a_2\otimes a_3)&=a_1\otimes a_2\otimes a_3+(-1)^{|a_1||a_2|}a_2\otimes a_1\otimes a_3+(-1)^{|a_1|(|a_2|+|a_3|)}a_2\otimes a_3\otimes a_1\cr
bat_{1,2}(a_2,a_1\otimes a_3)&=a_2\otimes a_1\otimes a_3+(-1)^{|a_2||a_1|}a_1\otimes a_2\otimes a_3+(-1)^{|a_2|(|a_1|+|a_3|)}a_1\otimes a_3\otimes a_2.
\endaligned
$$
Si $V=Vect\{a_1,a_2,a_3\}$, ces quatre tenseurs forment une base de l'espace image des battements de $Bat(1,2)$ et $Bat(2,1)$ dans $\otimes^3V$. Les tenseurs $a_1\otimes a_2\otimes a_3$ et $a_2\otimes a_1\otimes a_3$ forment une base d'un suppl\'ementaire, donc $\underline{a_1\otimes a_2\otimes a_3}$ et $\underline{a_2\otimes a_1\otimes a_3}$ forment une base de $\underline{\otimes^3V}$.\\

Ce qui remplacera l'alg\`ebre tensorielle $\oplus_{n\geq1}\otimes^nA[1]$ sera donc ici l'espace $\oplus_{n\geq1}\underline{\otimes^nR[1]}$. Le quotient tient compte bien s\^ur du nouveau degr\'e $\deg(a)=|a|-1$ des \'el\'ements de $R[1]$.\\

Cet espace est muni d'un cocrochet de Lie $\delta$ sur $\oplus_{n\geq1}\underline{\otimes^nR[1]}$ .On pose d'abord~:
$$\begin{aligned}
\delta(a_1\otimes\dots\otimes a_n)=~&\sum_{j=1}^{n-1}a_1\otimes\dots\otimes a_j\bigotimes a_{j+1}\otimes\dots\otimes a_n\\
&-\varepsilon_{\deg}\left(\begin{smallmatrix}1&\dots&n-j&n-j+1&\dots&n\\ j+1&\dots&n&1&\dots&j\end{smallmatrix}\right)a_{j+1}\otimes\dots\otimes a_n\bigotimes a_1\otimes\dots\otimes a_j.
\end{aligned}
$$

Cette formule permet de d\'efinir $\delta$ sur l'espace quotient $\underline{\otimes^n(R[1])}$ (voir \cite{[AAC3]}).\\

\begin{prop} {\rm (La structure de cog\`ebre)}

\

L'espace $\oplus_{n\geq1}\underline{\otimes^nR[1]}$ \'equipp\'e de $\delta$ est une cog\`ebre de Lie, c'est \`a dire que $\delta$ est coantisym\'etrique de degr\'e $0$ et v\'erifie l'identit\'e de coJacobi: si $\tau$ est la volte 
$$
\tau(a\otimes b)=(-1)^{\deg(a)\deg(b)}b\otimes a,
$$
alors
$$
\tau\circ\delta=-\delta,\qquad\Big(id^{\otimes3}+(\tau\otimes id)\circ(id\otimes\tau)+(id\otimes\tau)\circ(\tau\otimes id)\Big)\circ(\delta\otimes id)\circ\delta=0.
$$
\end{prop}

Cette cog\`ebre de Lie est libre. C'est \`a dire que si $(\mathcal{C},c)$ est une cog\`ebre de Lie nilpotente quelconque, tout $f:(\mathcal{C},c)\longrightarrow R[1]$ lin\'eaire se prolonge en $F:(\mathcal{C},c)\longrightarrow \oplus_{n\geq1}\underline{\otimes^nR[1]}$ qui est un morphisme de cog\`ebre. On en d\'eduit qu'on peut d\'efinir des cod\'erivations $D$ et des morphismes $F$ de cette structure \`a partir de leur `s\'erie de Taylor'.

Soit $F:\oplus_{n\geq1}\underline{\otimes^nR[1]}\longrightarrow\oplus_{n\geq1}\underline{\otimes^nS[1]}$ un morphisme de cog\`ebres de Lie~:
$$
(F\otimes F)\circ\delta=\delta\circ F,
$$
on suppose toujours $F$ homog\`ene de degr\'e 0. On appelle $F_n$ la projection sur $S[1]$ parall\`element \`a $\oplus_{n>1}\underline{\otimes^nS[1]}$ de la restriction de $F$ \`a $\underline{\otimes^n(R[1])}$~: $F_n$ est une application lin\'eaire de $\underline{\otimes^n(R[1])}$ dans $S[1]$.

De m\^eme soit $D:\oplus_{n\geq1}\underline{\otimes^nR[1]}\longrightarrow\oplus_{n\geq1}\underline{\otimes^nR[1]}$ une cod\'erivation de cog\`ebres de Lie~:
$$
(id\otimes D+D\otimes id)\circ\delta=\delta\circ D,
$$
on suppose $D$ homog\`ene de degr\'e $q$. On appelle $D_n$ la projection sur $R[1]$ parall\`element \`a $\oplus_{n>1}\underline{\otimes^nR[1]}$ de la restriction de $D$ \`a $\underline{\otimes^n(R[1])}$~: $D_n$ est une application lin\'eaire de $\underline{\otimes^n(R[1])}$ dans $R[1]$.\\

\begin{prop} {\rm (Reconstruction de $F$ et $D$)}

\

La suite d'applications $(F_n)$ (resp. $(D_n)$) permet de reconstruire $F$ (resp. $D$) de fa\c con unique. En posant $a_{[i,i+j]}=\underline{a_i\otimes a_{i+1}\otimes\dots\otimes a_{i+j}}$, on a explicitement~:
$$
F(a_{[1,n]})=\sum_{\begin{smallmatrix}k>0,~0<r_1,\dots,r_k\\ r_1+\dots+r_k=n\end{smallmatrix}}
\underline{F_{r_1}(a_{[1,r_1]})\otimes F_{r_2}(a_{[r_1+1,r_1+r_2]})\otimes\dots\otimes F_{r_k}(a_{[n-r_k+1,n]})}
$$
et
$$
D(a_{[1,n]})=\sum_{\begin{smallmatrix}0<r\\1\leq j\leq n-r\end{smallmatrix}}(-1)^{q\deg(a_{[1,j]})}\underline{a_{[1,j]}\otimes D_r(a_{[j+1,j+r]})\otimes a_{[j+r+1,n]}}.
$$
Plus pr\'ecis\'ement, toute suite d'applications $(\varphi_n)$ peut se relever d'une seule fa\c con en un morphisme (resp. une cod\'erivation).\\
\end{prop}

Voir \cite{[AAC3]} et \cite{[BGHHW]} pour la preuve de cette proposition.\\

La multiplication dans $R$ permet de transformer $\oplus_{n\geq1}\underline{\otimes^n(R[1])}$ en un complexe d'homologie. Posons $C_n(R)=\underline{\otimes^n(R[1])}$. La multiplication se remonte en une application $m~:\underline{\otimes^2R[1]}\longrightarrow R[1]$ comme ci-dessus~:
$$
m(\underline{a\otimes b})=(-1)^{\deg(a)}ab.
$$
Cette application est de degr\'e 1, anticommutative et antiassociative~:
$$\aligned
\deg(m(\underline{a\otimes b}))&=1+\deg(a)+\deg(b),\\ m(\underline{a\otimes b})&=-(-1)^{\deg(a)\deg(b)}m(\underline{b\otimes a}),\\ m(\underline{m(\underline{a\otimes b})\otimes c})&=-(-1)^{\deg(a)}m(\underline{a\otimes m(\underline{b\otimes c})}).\endaligned
$$
Elle se prolonge d'une fa\c con unique en une cod\'erivation de $(\oplus_{n\geq1}\underline{\otimes^nR[1]},\delta)$, toujours not\'ee $m$. On a~:
$$
m(\underline{a_0\otimes\dots\otimes a_n})=\sum_{j=0}^{n-1}(-1)^{\sum_{i<j}\deg(a_i)}\underline{a_0\otimes\dots\otimes m(a_j\otimes a_{j+1})\otimes\dots\otimes a_n}
$$
et donc~:
$$
(id\otimes m+m\otimes id)\circ\delta=\delta\circ m.
$$

La relation d'associativit\'e pour le produit est \'equivalente \`a
$$
m\circ m=0.
$$

On d\'efinit ainsi un complexe d'homologie appel\'e $C_\infty$ alg\`ebre enveloppante de $R$.\\

\begin{defn}{\rm(Alg\`ebre $C_\infty$ enveloppante)}

\

Soit $R$ une alg\`ebre gradu\'ee, associative et commutative. On appelle $C_\infty$ alg\`ebre enveloppante de $R$ la cog\`ebre de Lie $(\oplus_{n\geq1}\underline{\otimes^nR[1]},\delta)$ munie de la cod\'erivation $m$. On note cette alg\`ebre enveloppante
$$
C_\infty(R)=(\oplus_{n\geq1}\underline{\otimes^nR[1]},\delta,m).
$$
\end{defn}

La $C_\infty$ alg\`ebre enveloppante de $R$ est donc un complexe d'homologie~:
$$
0\stackrel{m}{\longleftarrow}C_0(R)\stackrel{m}{\longleftarrow}C_1(R)\stackrel{m}{\longleftarrow}\dots \stackrel{m}{\longleftarrow}C_{n-1}(R)\stackrel{m}{\longleftarrow}\dots
$$
o\`u $C_n(R)=\underline{\otimes^nR[1]}$
On note $Z_n(R)=\ker(m|_{C_n(R)})$ et $B_n(m)=\im(m|_{C_{n+1}(R)})$, on a donc les groupes d'homologie~:
$$
H_n(R)=Z_n(R)/B_n(R)
$$
pour tout $n$. On appelle cette homologie, l'homologie de Harrison de l'alg\`ebre commutative $R$.\\

Plus g\'en\'eralement,\\

Soit $R$ une alg\`ebre gradu\'ee, associative et commutative. Soit $M$ un $R$-bimodule gradu\'e et sym\'etrique~:
$$
v.a=(-1)^{|a||v|}a.v\qquad(a\in R,~~v\in M).
$$
On posera $m(\underline{a\otimes v})=(-1)^{\deg(a)}a.v$. L'homologie de Harrison du bimodule $M$ est d\'efinie de la fa\c con suivante.

On pose $C_n(R,M)=\underline{\sum_{j=0}^n\otimes^jR[1]\otimes M[1]\otimes\otimes^{n-j}R[1]}$, quotient de l'espace des $n$ cha\^\i nes de $R$ vue comme une alg\`ebre associative par l'espace des images de tous les battements $bat_{p,n-p}$ agissant sur le produit tensoriel avec un facteur $M[1]$. Le complexe d'homologie de Harrison de module $M$ est le complexe~:
$$
0\stackrel{\partial_0}{\longleftarrow}C_0(R,M)\stackrel{\partial_1}{\longleftarrow}C_1(R,M)\stackrel{\partial_2}{\longleftarrow}\dots \stackrel{\partial_{n-1}}{\longleftarrow}C_{n-1}(R,M)\stackrel{\partial_n}{\longleftarrow}\dots
$$
o\`u le bord $\partial_n$ est lin\'eaire, de degr\'e 1 et d\'efini comme le quotient par les images des battements du bord de Hochschild,c'est à dire par~:
$$\aligned
\partial_n(&\underline{a_0\otimes\dots\otimes a_{i-1}\otimes v_i\otimes a_{i+1}\otimes\dots\otimes a_n})=\cr
&\sum_{j=0}^{i-2}(-1)^{\sum_{k<j}\deg(a_k)}\underline{a_0\otimes\dots\otimes m(\underline{a_j\otimes a_{j+1}})\otimes\dots\otimes v_i\otimes\dots\otimes a_n}+\cr
&+(-1)^{\sum_{k<i-1}\deg(a_k)}\underline{a_0\otimes\dots\otimes m(\underline{a_{i-1}\otimes v_i})\otimes\dots\otimes a_n}\cr&+(-1)^{\sum_{k<i}\deg(a_k)}\underline{a_0\otimes\dots\otimes m(\underline{(v_i\otimes a_{i+1})}\otimes\dots\otimes a_n}+\cr
&+\sum_{j=i+1}^{n-1}(-1)^{\sum_{k<j}\deg(a_k)+\deg(v_i)}\underline{a_0\otimes\dots\otimes v_i\otimes\dots\otimes m(\underline{a_j\otimes a_{j+1}})\otimes\dots\otimes a_n}.
\endaligned
$$
Comme plus haut, le $n^{i\grave{e}me}$ groupe d'homologie du bimodule $M$ est
$$
H_n(R,M)=Z_n(R,M)/B_n(R,M),
$$
o\`u $Z_n(R,M)=\ker(\partial_n)$ et $B_n(R,M)=\im(\partial_{n+1})$.\\

De m\^eme la cohomologie de ce bimodule $M$ est d\'efinie ainsi. On pose
$$
C^n(R,M)=L(\underline{\otimes^nR[1]},M[1]),
$$
et on obtient le complexe de cohomologie de Harisson en posant~:
$$
0\stackrel{\partial^0}{\longrightarrow}C^0(R,M)\stackrel{\partial^1}{\longrightarrow}C^1(R,M)\stackrel{\partial^2}{\longrightarrow}\dots \stackrel{\partial^{n-1}}{\longrightarrow}C^{n-1}(R,M)\stackrel{\partial^n}{\longrightarrow}\dots
$$
avec, si $f\in C^n(R,M)$ est homog\`ene de degr\'e $\deg(f)$,
$$\aligned
(\partial^{n+1}f)(\underline{a_0\otimes\dots\otimes a_n})=&(-1)^{\deg(f)\deg(a_0)}m(a_0\underline\otimes f(\underline{a_1\otimes\dots\otimes a_n}))\cr&-\sum_{j=1}^{n-2}(-1)^{\sum_{k<j}\deg(a_k)}f(\underline{a_0\otimes\dots\otimes m(\underline{a_j\otimes a_{j+1}})\otimes\dots\otimes a_n})\cr
&+(-1)^{\sum_{k<n}\deg(a_k)}f(\underline{a_0\otimes\dots\otimes a_{n-1}})\underline\otimes a_n.
\endaligned
$$
Et la cohomologie de Harrison de $M$ est celle de ce complexe~:
$$
H^n(R,M)=Z^n(R,M)/B^n(R,M),\text{ o\`u } Z^n(R,M)=\ker(\partial^{n+1})\text{ et } B^n(A,M)=\im(\partial^n).
$$

Remarquons que si $R$ et $M$ ne sont pas gradu\'es, on retrouve l'homologie et la cohomologie de Harrison usuelles du $R$ bimodule sym\'etrique $M$.\\

\begin{prop}{\rm(Homologie et cohomologie de Harrison et $C_\infty$ alg\`ebres)}

\

Soit $R$ une alg\`ebre associative, commutative et gradu\'ee. Soit $M$ un $R$-bimodule sym\'etrique gradu\'e. On d\'efinit une nouvelle alg\`ebre associative et commutative, not\'ee $S=R\ltimes M$, produit semi-direct de $R$ par $M$ en munissant $R\oplus M$ de la multiplication~:
$$
(a+u)(b+v)=(ab+(ub+av))\qquad (a,~b\in R,~~u,~v\in M).
$$
Alors
\begin{itemize}
\item[1.] Le complexe d'homologie de Harrison du bimodule sym\'etrique $M$ est un sous complexe de $(C_\infty(R\ltimes M),m_{R\ltimes M})$.\\

\item[2.] Un morphisme de cog\`ebres $F:(\oplus_{n\geq1}\underline{\otimes^nR[1]},\delta_R)\longrightarrow(\oplus_{n\geq1}\underline{\otimes^nS[1]},\Delta_S)$, de degr\'e $0$ d\'efini par $F_1=\iota+c_1$, $F_n=c_n$ o\`u $\iota$ est l'injection canonique de $R$ dans $R\ltimes M$ et $c_j\in L(\otimes^jR[1],M[1])$ est un morphisme
$$
F~:(C_\infty(R),\delta_R,m_R)\longrightarrow(C_\infty(S),\delta_S,m_S)
$$
de $C_\infty$ alg\`ebres si et seulement si $m_S\circ F=F\circ m_R$, si et seulement si $\partial c_j=0$ ($j\geq1$) o\`u $\partial$ est le cobord de Harrison.

\item[3.] Un tel morphisme est dit trivial s'il existe une suite d'applications $b=(b_n)$, de degr\'e -1 telle que $c=m_B\circ(\iota\otimes b-b\otimes\iota)-b\circ m_A$, ce qui est \'equivalent à $c_1=0$ et $c_n=\partial b_{n-1}$, pour tout $n>1$.\\
\end{itemize}
\end{prop}

\

\noindent
{\bf Preuve}

\

\noindent
1. Comme pour les $A_\infty $ alg\'ebre, $S=R\oplus M$, alors,
$$
C_n(R,M)=\underline{\sum_{j=0}^{n-1}\otimes^jR\otimes M\otimes^{n-1-j}R}\subset\underline{\otimes^{n-1}S}=C_n(S) 
$$
et ${m_S}_{|C_n(R,M)}=\partial_n$.

D'o\`u le complexe d'homologie de Harrison du bimodule $M$ est un sous complexe de $\Big(C_\infty(R\ltimes M), m_{R\ltimes M}\Big)$.\\

\noindent
2. Soit $F:\Big(\displaystyle\bigoplus_{n\geq1}\underline{\otimes^n R[1]},\Delta_R,m_R\Big)\longrightarrow\Big(\displaystyle\bigoplus_{n\geq1}
\underline{\otimes^n S[1]},\Delta_S,m_S\Big)$ un morphisme de $C_\infty$-alg\`ebres tel que $F_1=\iota+c_1$ et $F_n=c_n$. Alors, on a
$$
F(\underline{a_1\otimes\dots \otimes a_n})=\sum_{\begin{smallmatrix}k>0,~0<r_1,\dots,r_k\\ r_1+\dots+r_k=n\end{smallmatrix}}\underline{F_{r_1}(a_1\otimes\dots \otimes a_{r_1})
\otimes \dots\otimes F_{r_k}(a_{n-r_k+1}\otimes\dots \otimes a_n)}.
$$
et $m_S\circ F=F\circ m_R$. On va d\'emontrer que $\partial c_j=0$ par r\'ecurrence sur $j$.

Pour $j=1$, $F(\underline{a\otimes b})=\underline{F_1(a)\otimes F_1(b)}+ F_2(\underline{a\otimes b})$.
Donc,
$$\aligned
m_S\circ F(\underline{a\otimes b})&=m_S(\underline{F_1(a)\otimes F_1(b)})=(-1)^{dega}(ab+c_1(a)b+ac_1(b))\\
&=F\circ m_R(\underline{a\otimes b})= (-1)^{dega}(ab+c_1(ab)).
\endaligned
$$
D'o\`u $\partial c_1(\underline{a\otimes b})=0$.\\

Supposons que $\partial c_j=0$, pour tout $ j\leq n-2$ et montrons que $\partial c_{n-1}=0$ ($n\geq 2$). On a
\begin{align*}
F\circ m_R&(\underline{a_0\otimes\dots\otimes a_n})=F(\sum_{j=0}^{n-1}(-1)^{\sum_{i<j}\deg(a_i)}\underline{a_0\otimes\dots\otimes m_R(\underline{a_j\otimes a_{j+1}})\otimes\dots\otimes a_n})\\
=~&\sum_{\begin{smallmatrix}k>0,~0<r_1,\dots,r_k\\ r_1+\dots+r_k=n\end{smallmatrix}}\sum_{j=0}^{n-1}(-1)^{\sum_{i<r_1+...+r_{ j-1}+1}\deg(a_i)}\underline{F_{r_1}(\underline{a_1\otimes\dots \otimes a_{r_1}})
\otimes }\\&\underline{\otimes\dots\otimes F_{r_j}\circ m_R(\underline{a_{r_1+...+r_{ j-1}+1}\otimes\dots \otimes a_{r_j}})
\otimes \dots\otimes F_{r_k}(\underline{a_{n-r_k+1}\otimes\dots \otimes a_n})}
\end{align*}
et
\begin{align*}
m_S \circ F&(\underline{a_0\otimes\dots\otimes a_n})=\\
=&\sum_{\begin{smallmatrix}k>0,~0<r_1,\dots,r_k\\ r_1+\dots+r_k=n\end{smallmatrix}}\sum_{j=0}^{n-1}(-1)^{\sum_{i<r_1+...+r_{j-1}+1}\deg(a_i)} \underline{F_{r_1}(\underline{a_1\otimes\dots \otimes a_{r_1}})
\otimes\dots\otimes}\\
&\underline{\otimes m_S\Big( F_{r_j}(\underline{a_{r_1+...+r_{ j-1}+1}\otimes\dots \otimes a_{r_1+...+r_j}})\underline{\otimes} F_{r_{j+1}}(\underline{a_{r_1+...+r_{ j}+1}\otimes\dots \otimes a_{r_1+...+r_{j+1}}})\Big)}\\&
\underline{\otimes \dots\otimes F_{r_k}(\underline{a_{n-r_k+1}\otimes\dots \otimes a_n})}.
\end{align*}

Dans la somme pr\'ec\'edente, si $r_j>1$ et $r_{j+1}>1$ on a
$$
m_S\Big( F_{r_j}(\underline{a_{r_1+...+r_{ j-1}+1}\otimes\dots \otimes a_{r_1+...+r_j}})\otimes F_{r_{j+1}}(\underline{a_{r_1+...+r_{ j}+1}\otimes\dots \otimes a_{r_1+...+r_{j+1}}})\Big)=0.
$$
De plus, on a si $r_j=1$,
\begin{align*}
m_S\Big(&\underline{F_1(a_{r_1+...+r_{ j-1}+1})\otimes F_{r_{j+1}}(\underline{a_{r_1+...+r_{ j}+1}\otimes\dots \otimes a_{r_1+...+r_{j+1}}})}\Big)=\\
&\hskip 4cm=m_S\Big(\underline{a_{r_1+...+r_{ j-1}+1}\otimes F_{r_{j+1}}(\underline{a_{r_1+...+r_{ j}+1}\otimes\dots \otimes a_{r_1+...+r_{j+1}}})}\Big)
\end{align*}
et si $r_{j+1}=1$,
\begin{align*}
m_S\Big(&\underline{F_{r_j}(\underline{a_{r_1+...+r_{j-1}+1}\otimes\dots\otimes a_{r_1+...+r_j}})\underline{\otimes} F_1(a_{r_1+...+r_{j}+1})}\Big)\\
&\hskip 4cm=m_S\Big(\underline{F_{r_j}(\underline{a_{r_1+...+r_{ j-1}+1}\otimes\dots \otimes a_{r_1+...+r_j}})\otimes a_{r_1+...+r_{ j}+1}}\Big).
\end{align*}

Dans l'expression $(m_S \circ F-F\circ m_R)(\underline{a_0\otimes\dots\otimes a_n})$, les termes o\`u $r_j\leq n-2$ et les termes o\`u $r_{j+1}\leq n-2$ disparaissent gr\^ace à notre hypoth\`ese de r\'ecurrence $\partial c_{r_j}=0$, pour tout $r_j\leq n-2$. Il ne reste que les termes o\`u $r_j=1$ et $r_{j+1}=n-1$ ou $r_j=n-1$ et $r_{j+1}=1$. Donc
$$
\begin{aligned}
(m_S \circ F&-F\circ m_R)(\underline{a_0\otimes\dots\otimes a_n})=m_S\big(c_{n-1}(\underline{a_0\otimes\dots\otimes a_{n-1}})\underline{\otimes} a_n\big)\\
&+m_S\big(a_0\underline{\otimes} c_{n-1}(\underline{a_1\otimes\dots\otimes a_{n}})\big)\cr
&-\sum_{j=1}^{n-1}(-1)^{\sum_{i<j}deg a_i}c_{n-1}\big(\underline{a_0\otimes\dots\otimes a_{j-1}\otimes m_R
(a_j\underline{\otimes} a_{j+1})\otimes a_{j+2}\otimes\dots\otimes a_n}\big)\cr
=~&(-1)^{\sum_{i\leq n-1}deg a_i}c_{n-1}(\underline{a_0\otimes\dots\otimes a_{n-1}}).a_n+(-1)^{deg a_0}a_0. c_{n-1}(\underline{a_1\otimes\dots\otimes a_{n}})\cr
&-\sum_{j=1}^{n-1}(-1)^{\sum_{i\leq j}deg a_i}c_{n-1}\big(\underline{a_0\otimes\dots\otimes a_{j-1}\otimes(a_j.a_{j+1})\otimes a_{j+2}\otimes\dots \otimes a_n}\big)\cr
=~&(\partial c_{n-1})(\underline{a_0\otimes\dots\otimes a_n})=0.
\end{aligned}
$$
On retrouve l'op\'erateur de cobord de Harrison $\partial $.

\

\noindent
3. Si $c=m_S\circ(\iota\otimes b-b\otimes\iota)-b\circ m_A$, alors,
$$
\begin{aligned}
c_{n-1}(&\underline{a_0\otimes\dots\otimes a_n})=\\
=~&m_S(a_0\underline{\otimes} b_{n-2}(\underline{a_1\otimes\dots\otimes a_n})+b_{n-2}(\underline{a_0\otimes\dots a_{n-1}})\underline{\otimes} a_n)\\
&-\sum_{j=0}^{n-1}(-1)^{\sum_{i\leq j}deg a_i}b_{n-2}(\underline{a_0\otimes\dots \otimes a_ja_{j+1}\otimes \dots\otimes a_n})\\
=~&(-1)^{deg a_0}a_0. b_{n-2}(\underline{a_1\otimes\dots\otimes a_n})+(-1)^{\sum_{i\leq n-1}deg a_i}b_{n-2}(\underline{a_0\otimes\dots \otimes a_{n-1}}).a_n\\
&-\sum_{j=0}^{n-1}(-1)^{\sum_{i\leq j}deg a_i}b_{n-2}(\underline{a_0\otimes\dots \otimes a_ja_{j+1}\otimes \dots\otimes a_n})\\
=~&\partial b_{n-2}(\underline{a_0\otimes\dots\otimes a_n}).
\end{aligned}
$$


\section{Alg\`ebre $L_\infty$ enveloppante d'une alg\`ebre de Lie}\label{sec4}\

Dans ce paragraphe, $\mathfrak g$ est une alg\`ebre de Lie gradu\'ee sur $\R$. La relation d'antisym\'etrie du crochet~:
$Antisym(a_1\otimes a_2)=[a_1,a_2]+(-1)^{|a_1||a_2|}[a_2,a_1]=0$ est bien s\^ur sym\'etrique, de plus, la relation de Jacobi $Jac(a_1,a_2,a_3)=0$ o\`u~:
$$
Jac(a_1,a_2,a_3)=(-1)^{|a_1||a_3|}\big[[a_1,a_2],a_3\big]+(-1)^{|a_2||a_1|}\big[[a_2,a_3],a_1\big]+(-1)^{|a_3||a_2|}\big[[a_3,a_1],a_2\big]
$$
v\'erifie~:
$$
\dim Vect\{Jac\circ\sigma,~~\sigma\in\mathfrak{S}_3\}=1.
$$
En fait $Jac\circ\sigma(a_1,a_2,a_3)=0$ pour tout $\sigma$ est \'equivalent, en tenant compte de l'antisym\'etrie, \`a la relation $Jac(a_1,a_2,a_3)=0$.\\

Apr\`es d\'ecalage des degr\'es de 1, on doit donc consid\'erer des tenseurs compl\`etements sym\'etriques pour le degr\'e $\deg$.

Ce qui remplacera l'alg\`ebre tensorielle $\oplus_{n\geq1}\otimes^nA[1]$ sera donc ici l'espace $\oplus_{n\geq1}S^n(\mathfrak g[1])$. La sym\'etrie tient compte bien s\^ur du nouveau degr\'e $\deg(a)=|a|-1$ des \'el\'ements de $\mathfrak g[1]$. Un \'el\'ement de $S^n(\mathfrak g[1])$ est not\'e $a_1\pt a_2\pt\dots\pt a_n$.\\

Cet espace est muni d'une comultiplication $\Delta$ d\'efinie par~:
$$
\Delta(a_1\pt\dots\pt a_n)=
\sum_{\begin{smallmatrix}
I\sqcup J=\{1,\dots n\}\\
|I|>0,|J|>0
\end{smallmatrix}}
 \varepsilon_{\deg(a)}
 \left(\begin{smallmatrix}
  \{1,\dots n\}\\
   I,J
   \end{smallmatrix}\right)a_I\otimes a_J,
$$
o\`u $a_I=a_{i_1}\pt a_{i_2}\pt\dots\pt a_{i_{|I|}}$ si $I=\{i_1<i_2<\dots<i_{|I|}\}$. $\Delta$ est ainsi bien d\'efini sur $S^n(\mathfrak g[1])$ (voir \cite{[AMM]}).\\

\begin{prop} {\rm (La structure de cog\`ebre)}

\

L'espace $\oplus_{n\geq1}S^n(\mathfrak g[1])$ \'equipp\'e de $\Delta$ est une cog\`ebre cocommutative et coassociative, c'est \`a dire que $\Delta$ est de degr\'e $0$ et v\'erifie, si $\tau$ est la volte 
$$
\tau(a\otimes b)=(-1)^{\deg(a)\deg(b)}b\otimes a,
$$
alors
$$
\tau\circ\Delta=\Delta,\qquad(id\otimes\Delta)\circ\Delta=(\Delta\otimes id)\circ\Delta.
$$
\end{prop}

Cette cog\`ebre est libre, en particulier, on peut d\'efinir des cod\'erivations $Q$ et des morphismes $F$ de cette structure \`a partir de leurs `s\'erie de Taylor'.

Si $\mathfrak h$ est une alg\'ebre de lie gradu\'ee, soit $F~:~\oplus_{n\geq1}S^n(\mathfrak g[1])\longrightarrow\oplus_{n\geq1}S^n(\mathfrak h[1])$ un morphisme de cog\`ebres~:
$$
(F\otimes F)\circ\Delta=\Delta\circ F.
$$
On suppose toujours $F$ homog\`ene de degr\'e 0. On appelle $F_n$ la projection sur $\mathfrak h[1]$ parall\`element \`a $\oplus_{n>1}S^n(\mathfrak h[1])$ de la restriction de $F$ \`a $S^n(\mathfrak g[1])$~: $F_n$ est une application lin\'eaire de $S^n(\mathfrak g[1])$ dans $\mathfrak h[1]$.

De m\^eme soit $D:\oplus_{n\geq1}S^n(\mathfrak g[1])\longrightarrow\oplus_{n\geq1}S^n(\mathfrak g[1])$ une cod\'erivation de $\Delta$~:
$$
(id\otimes D+D\otimes id)\circ\Delta=\Delta\circ D.
$$
On suppose $D$ homog\`ene de degr\'e $q$. On appelle $D_n$ la projection sur $\mathfrak g[1]$ parall\`element \`a $\oplus_{n>1}S^n(\mathfrak g[1])$ de la restriction de $D$ \`a $S^n(\mathfrak g[1])$~: $D_n$ est une application lin\'eaire de $S^n(\mathfrak g[1])$ dans $\mathfrak g[1]$.\\

\begin{prop} {\rm (Reconstruction de $F$ et $D$)}

\

La suite d'applications $(F_n)$ (resp. $(D_n)$) permet de reconstruire $F$ (resp. $D$) de fa\c con unique. On a explicitement~:
$$
F(a_1\pt\dots\pt a_n)=\sum_{j>0}\frac{1}{j!}\sum_{\begin{smallmatrix}I_1\sqcup\cdots\sqcup I_j=\{1,\ldots,n\}\\ I_1\dots I_j\neq\emptyset\end{smallmatrix}}\varepsilon_{\deg(a)}\left(\begin{smallmatrix}1,\dots,n\\ I_1,\dots,I_j\end{smallmatrix}\right) F_{|I_1|}(a_{I_1})\pt\dots\pt F_{|I_j|}(a_{I_j}).
$$
et
$$
D(a_1\pt\dots\pt a_n)=\sum_{\begin{smallmatrix}I\sqcup J=\{1,\dots,n\}\\ I,J\neq\emptyset\end{smallmatrix}} \varepsilon_{\deg(a)}\left(\begin{smallmatrix}1,\dots,n\\ I, J\end{smallmatrix}\right)\bigl(D_{|I|}(a_I)\bigr)\pt a_J.
$$
Plus pr\'ecis\'ement, toute suite d'applications $(\varphi_n)$ peut se relever d'une seule fa\c con en un morphisme (resp. une cod\'erivation).\\
\end{prop}

Voir \cite{[AMM]} et \cite{[AAC3]} pour la preuve de cette proposition.\\

Le crochet de Lie dans $\mathfrak g$ permet de transformer $\oplus_{n\geq1}S^n(\mathfrak g[1])$ en un complexe d'homologie. Posons $C_n(\mathfrak g)=S^n(\mathfrak g[1])$. Le crochet de Lie se remonte en une application $\ell~:S^2(\mathfrak g[1])\longrightarrow \mathfrak g[1]$ comme ci-dessus~:
$$
\ell(a\pt b)=(-1)^{\deg(a)}[a,b].
$$
Cette application est de degr\'e 1, commutative et v\'erifie Jacobi~:
$$\aligned
&\deg(\ell(a\pt b))=1+\deg(a)+\deg(b),\\ &\ell(a\pt b)=(-1)^{\deg(a)\deg(b)}\ell(b\pt a),\\ &(-1)^{\deg(a_1)\deg(a_3)}\ell\big(\ell(a_1\pt a_2)\pt a_3\big)+(-1)^{\deg(a_2)\deg(a_1)}\ell\big(\ell(a_2\pt a_3)\pt a_1\big)+\\ &\hskip 2cm+(-1)^{\deg(a_3)\deg(a_2)}\ell\big(\ell(a_3\pt a_1)\pt a_2\big)=0.\endaligned
$$
Elle se prolonge d'une fa\c con unique en une cod\'erivation de $(\oplus_{n\geq1}S^n(\mathfrak g[1]),\Delta)$, toujours not\'ee $\ell$. On a~:
$$
\ell(a_1\pt\dots\pt a_n)=\sum_{i<j}\varepsilon_{\deg(a)}\left(\begin{smallmatrix}1\ \dots\ n\\ i,j,1\dots\hat{\imath}\hat{\jmath}\dots n\end{smallmatrix}\right)\ell(a_i\pt a_j)\pt a_1\pt\dots\hat\imath\dots\hat\jmath\dots\pt a_n
$$
et donc~:
$$
(id\otimes\ell+\ell\otimes id)\circ\Delta=\Delta\circ\ell.
$$

La relation de Jacobi pour le crochet est \'equivalente \`a
$$
\ell\circ\ell=0.
$$

On d\'efinit ainsi un complexe d'homologie appel\'e $L_\infty$ alg\`ebre enveloppante de $\mathfrak g$.\\

\begin{defn}{\rm(Alg\`ebre $L_\infty$ enveloppante)}

\

Soit $\mathfrak g$ une alg\`ebre de Lie gradu\'ee. On appelle $L_\infty$ alg\`ebre enveloppante de $\mathfrak g$ la cog\`ebre coassociative et cocommutative $(\oplus_{n\geq1}S^n(\mathfrak g[1]),\Delta)$ munie de la cod\'erivation $\ell$. On note cette alg\`ebre enveloppante
$$
L_\infty(\mathfrak g)=(\oplus_{n\geq1}S^n(\mathfrak g[1]),\Delta,\ell).
$$
\end{defn}

La $L_\infty$ alg\`ebre enveloppante de $\mathfrak g$ est donc un complexe d'homologie~:
$$
0\stackrel{\ell}{\longleftarrow}C_1(\mathfrak g)\stackrel{\ell}{\longleftarrow}C_2(\mathfrak g)\stackrel{\ell}{\longleftarrow}\dots \stackrel{\ell}{\longleftarrow}C_{n-1}(\mathfrak g)\stackrel{\ell}{\longleftarrow}\dots
$$
o\`u $C_n(\mathfrak g)=S^n(\mathfrak g[1])$. On note $Z_n(\mathfrak g)=\ker(\ell|_{C_n(\mathfrak g)})$ et $B_n(\mathfrak g)=\im(\ell|_{C_{n+1}(\mathfrak g)})$, on a donc les groupes d'homologie~:
$$
H_n(\mathfrak g)=Z_n(\mathfrak g)/B_n(\mathfrak g)
$$
pour tout $n$. Il est facile de voir que cette homologie est l'homologie de Chevalley de l'alg\`ebre de Lie $\mathfrak g$.\\

Plus g\'en\'eralement,\\

Soit $\mathfrak g$ une alg\`ebre de Lie gradu\'ee. Soit $M$ un $\mathfrak g$-module gradu\'e. On posera
$$
\ell(a\pt v)=(-1)^{\deg(a)}av=(-1)^{\deg(a)\deg(v)}\ell(v\pt a).
$$

L'homologie de Chevalley du module $M$ est d\'efinie de la fa\c con suivante.

On pose $C_n(\mathfrak g,M)=M[1]\pt S^n(\mathfrak g[1])$, c'est un sous espace de $S^{n+1}(\mathfrak g[1]\oplus M[1])$.

Le complexe d'homologie de Chevalley est le complexe~:
$$
0\stackrel{\partial_0}{\longleftarrow}C_0(\mathfrak g,M)\stackrel{\partial_1}{\longleftarrow}C_1(\mathfrak g,M)\stackrel{\partial_2}{\longleftarrow}\dots \stackrel{\partial_{n-1}}{\longleftarrow}C_{n-1}(\mathfrak g,M)\stackrel{\partial_n}{\longleftarrow}\dots
$$
o\`u le bord $\partial_n$ est lin\'eaire, de degr\'e 1 et d\'efini, en posant $v\pt a_1\pt\dots\pt a_n=x_0\pt\dots\pt x_n$, par~:
$$
\partial_n(x_0\pt\dots\pt x_n)=\sum_{i<j}\varepsilon_{\deg(x)}
\left(\begin{smallmatrix}
0\ \dots \ n \\
 i,j,0\dots\hat{\imath}\hat{\jmath}\dots n
\end{smallmatrix}\right)\ell(x_i\pt x_j)\pt x_0\pt\dots\hat{\imath}\dots\hat{\jmath}\dots\pt x_n
$$

Comme plus haut, le $n^{i\grave{e}me}$ groupe d'homologie du module $M$ est
$$
H_n(\mathfrak g,M)=Z_n(\mathfrak g,M)/B_n(\mathfrak g,M),
$$
o\`u $Z_n(\mathfrak g,M)=\ker(\partial_n)$ et $B_n(\mathfrak g,M)=\im(\partial_{n+1})$.\\

De m\^eme la cohomologie de ce module $M$ est d\'efinie ainsi. On pose
$$
C^n(\mathfrak g,M)=L(S^n(\mathfrak g[1]),M[1]),
$$
et on obtient le complexe de cohomologie de Chevalley en posant~:
$$
0\stackrel{\partial^0}{\longrightarrow}C^0(\mathfrak g,M)\stackrel{\partial^1}{\longrightarrow}C^1(\mathfrak g,M)\stackrel{\partial^2}{\longrightarrow}\dots \stackrel{\partial^{n-1}}{\longrightarrow}C^{n-1}(\mathfrak g,M)\stackrel{\partial^n}{\longrightarrow}\dots
$$
avec, si $f\in C^n(\mathfrak g,M)$ est homog\`ene de degr\'e $\deg(f)$~:
$$
\aligned
\partial_nf(a_1\pt\dots\pt a_{n+1})=~&\sum_i(-1)^{\deg(a_i)\deg(f)}\varepsilon_{\deg(a)}\left(\begin{smallmatrix}1~~~\dots~~~ n+1\\ i,1,\dots\hat\imath\dots n+1\end{smallmatrix}\right)\ell(a_i\pt f(a_1\pt\dots\hat{\imath}\dots\pt a_{n+1}))\\
&-\sum_{i<j}\varepsilon_{\deg(a)}\left(\begin{smallmatrix}1,~~~\dots~~~ ,n+1\\ i,j,1\dots\hat\imath\hat\jmath\dots n+1\end{smallmatrix}\right) f(\ell(a_i\pt a_j)\pt a_1\dots\hat\imath\dots\hat\jmath\dots\pt a_{n+1}).
\endaligned
$$

La cohomologie de Chevalley de $M$ est celle de ce complexe~:
$$
H^n(\mathfrak g,M)=Z^n(\mathfrak g,M)/B^n(\mathfrak g,M),\text{ o\`u } Z^n(\mathfrak g,M)=\ker(\partial^{n+1})\text{ et } B^n(\mathfrak g,M)=\im(\partial^n).
$$

Remarquons que si $\mathfrak g$ et $M$ ne sont pas gradu\'es $deg(f)=0,deg(a)=deg(m)=-1$, on retrouve l'homologie et la cohomologie de Chevalley usuelles du $\mathfrak g$ module $M$.\\

\begin{prop}{\rm(Homologie et cohomologie de Chevalley et $L_\infty$ alg\`ebres)}

\

Soit $\mathfrak g$ une alg\`ebre de Lie gradu\'ee. Soit $M$ un $\mathfrak g$-module gradu\'e. On d\'efinit une nouvelle alg\`ebre de Lie gradu\'ee not\'ee $\mathfrak h=\mathfrak g\ltimes M$, produit semi-direct de $\mathfrak g$ par $M$ en munissant $\mathfrak g\oplus M$ du crochet~:
$$
[(a+u),(b+v)]=([a,b]+(av-(-1)^{|b||u|}bu))\qquad (a,~b\in\mathfrak g,~~u,~v\in M).
$$
Alors
\begin{itemize}
\item[1.] Le complexe d'homologie de Chevalley du module $M$ est un sous complexe de $(L_\infty(\mathfrak g\ltimes M),\ell_{\mathfrak g\ltimes M})$.\\

\item[2.] Un morphisme de cog\`ebres $F:(\oplus_{n\geq1}S^n(\mathfrak g[1]),\Delta_{\mathfrak g})\longrightarrow(\oplus_{n\geq1}S^n(\mathfrak h[1]),\Delta_{\mathfrak h})$, de degr\'e $0$ d\'efini par $F_1=\iota+c_1$, $F_n=c_n$ o\`u $\iota$ est l'injection canonique de $\mathfrak g$ dans $\mathfrak g\ltimes M$ et $c_j\in L(S^j(\mathfrak g[1]),M[1])$ est un morphisme
$$
F~:~(L_\infty(\mathfrak g),\Delta_{\mathfrak g},\ell_{\mathfrak g})\longrightarrow(L_\infty(\mathfrak h),\Delta_{\mathfrak h},\ell_{\mathfrak h})
$$
de $L_\infty$ alg\`ebres si et seulement si $\ell_{\mathfrak h}\circ F=F\circ\ell_{\mathfrak g}$, si et seulement si $\partial c_j=0$ ($j\geq1$) o\`u $\partial$ est le cobord de Chevalley.

\item[3.] Un tel morphisme est dit trivial s'il existe une suite d'applications $b=(b_n)$, de degr\'e -1 telle que $c=\ell_{\mathfrak h}\circ(\iota\bullet b)-b\circ\ell_{\mathfrak g}$, o\`u par d\'efinition
$$
(\iota\bullet b)(a_1.\cdots.a_n)=\sum_i(-1)^{\deg(a_i)}\varepsilon_{\deg(a)}\left(\begin{smallmatrix}1~~~\dots~~~ n\\ i,1,\dots\hat\imath\dots n\end{smallmatrix}\right) a_i\pt b(a_1\pt\dots\hat{\imath}\dots\pt a_n)
$$
(qui est une g\'en\'eralisation de la formule d\'efinissant un morphisme en tenant compte du degr\'e de $b$). Ceci est \'equivalent à $c_1=0$ et $c_n=\partial b_{n-1}$, pour tout $n>1$.\\
\end{itemize}
\end{prop}

\

\noindent
{\bf Preuve}

\

\noindent
1. Par un calcul direct, on a~:
$$
C_n(\mathfrak g, M)=M[1].S^n(\mathfrak g[1])\subset S^{n+1}(\mathfrak g[1]\oplus M[1])~~\text{ et }~~\ell_{\mathfrak g\ltimes M}/C_n(\mathfrak g, M)=\partial_n.
$$
D'o\`u le complexe d'homologie de Chevalley du module $M$ est un sous complexe de $(L_\infty(\mathfrak g\ltimes M),\ell_{\mathfrak g\ltimes M})$.

\

\noindent
2. Soit $F:(\oplus_{n\geq1}S^n(\mathfrak g[1]),\Delta_{\mathfrak g}, \ell_{\mathfrak g})\longrightarrow(\oplus_{n\geq1}S^n(\mathfrak h[1]),\Delta_{\mathfrak h},\ell_{\mathfrak h})$ un morphisme de $L_\infty$ alg\`ebres de degr\'e $0$ tel que $F_1=\iota+c_1$ et $F_n=c_n$. Montrons par r\'ecurrence sur $j$ que $\partial c_j=0$.

On a $\ell_{\mathfrak h}\circ F=F\circ\ell_{\mathfrak g}$. D'une part,
\begin{align*}
\ell_{\mathfrak h}\circ F(a\otimes b)&=\ell_{\mathfrak h}(F_1(a).F_1(b)+F_2(a\otimes b))\\
&=\ell_{\mathfrak h}(F_1(a).F_1(b)),
\end{align*}
d'autre part :
$$
F\circ \ell_{\mathfrak g}(a\otimes b)=\ell_{\mathfrak g}(a\otimes b)+c_1(\ell_{\mathfrak g}(a\otimes b)).
$$
D'o\`u $\partial c_1(a\otimes b)=0$.

Supposons le r\'esultat vrai \`a l'ordre $n-2$ et montrons qu'il est vrai \`a l'ordre $n-1$. On a
\begin{align*}
&\ell_{\mathfrak h}\circ F(a_1\otimes\dots\otimes a_n)=\\
=&\ell_{\mathfrak h}\Big(\sum_{j>0}\frac{1}{j!}\sum_{\begin{smallmatrix}I_1\sqcup\cdots\sqcup I_j=\{1,\ldots,n\}\\ I_1\dots I_j\neq\emptyset\end{smallmatrix}}\varepsilon_{\deg(a)}\left(\begin{smallmatrix}1,\dots,n\\ I_1,\dots,I_j\end{smallmatrix}\right) F_{|I_1|}(a_{I_1})\pt\dots\pt F_{|I_j|}(a_{I_j})\Big)\\
=&\sum_{j>0}\frac{1}{j!}\hskip-0.3cm\sum_{\begin{smallmatrix}I_1\sqcup\cdots\sqcup I_j=\{1,\ldots,n\}\\ I_1\dots I_j\neq\emptyset\end{smallmatrix}}\sum_{s<k} \varepsilon_{\deg(a)}\left(\begin{smallmatrix}1,\dots,n\\ I_s,I_k,I_1,\dots \hat{sk}\dots,I_j\end{smallmatrix}\right)\ell_{\mathfrak h}\Big(F_{|I_s|}(a_{I_s}),F_{|I_k|}(a_{I_k})\Big)F_{|I_1|}(a_{I_1})\pt\dots\pt F_{|I_j|}(a_{I_j}).
\end{align*}

Dans cette somme, si $|I_s|>1$ et $|I_k|>1$, alors $\ell_{\mathfrak h}\Big(F_{|I_s|}(a_{I_s}),F_{|I_k|}(a_{I_k})\Big)=0$. Il reste donc :
\begin{align*}
&\sum_{j>0}\frac{1}{j!}\hskip-0.3cm\sum_{\begin{smallmatrix}I_1\sqcup\cdots\sqcup I_j=\{1,\ldots,n\}\\ I_1\dots I_j\neq\emptyset\end{smallmatrix}}\sum_{s\neq k\atop |I_s|=1} \varepsilon_{\deg(a)}\left(\begin{smallmatrix}1,\dots,n\\ s,I_k,I_1,\dots \hat{sk}\dots,I_j\end{smallmatrix}\right)\ell_{\mathfrak h}\Big(a_s,F_{|I_k|}(a_{I_k})\Big)F_{|I_1|}(a_{I_1})\pt\dots\pt F_{|I_j|}(a_{I_j})=\\
=&\sum_{j>0}\frac{1}{j!}\hskip-0.3cm\sum_{\begin{smallmatrix}I_1\sqcup\cdots\sqcup I_j=\{1,\ldots,n\}\\ I_1\dots I_j\neq\emptyset\end{smallmatrix}}\sum_{s\neq k\atop |I_s|=1, |I_k|\leq n-2} \hskip-0.5cm\varepsilon_{\deg(a)}\left(\begin{smallmatrix}1,\dots,n\\ s,I_k,I_1,\dots \hat{sk}\dots,I_j\end{smallmatrix}\right)\ell_{\mathfrak h}\Big(a_s,F_{|I_k|}(a_{I_k})\Big)F_{|I_1|}(a_{I_1})\pt\dots\pt F_{|I_j|}(a_{I_j})\\
&\hskip 1cm+\sum_{j>0}\frac{1}{j!}\hskip-0.3cm\sum_{\begin{smallmatrix}I_1\sqcup\cdots\sqcup I_j=\{1,\ldots,n\}\\ I_1\dots I_j\neq\emptyset\end{smallmatrix}}\sum_{s} \varepsilon_{\deg(a)}\left(\begin{smallmatrix}1,\dots,n\\ s,1,\dots \hat{s}\dots,n\end{smallmatrix}\right)\ell_{\mathfrak h}\Big(a_s,F_{n-1}(a_1\otimes\dots\hat{s}\dots\otimes a_{n})\Big).
\end{align*}
Car $\ell_{\mathfrak h}\Big(F_1(a_s),F_{n-1}(a_1\otimes\dots\hat{s}\dots\otimes a_{n})\Big)=\ell_{\mathfrak h}\Big(a_s,F_{n-1}(a_1\otimes\dots\hat{s}\dots\otimes a_{n})\Big)$.

D'autre part, on a
$$
F\circ\ell_{\mathfrak g}(a_1\otimes\dots\otimes a_n)=F\Big(\sum_{i<j}\varepsilon_{\deg(a)}\left(\begin{smallmatrix}1\ \dots\ n\\ i,j,1\dots\hat{\imath}\hat{\jmath}\dots n\end{smallmatrix}\right)\ell_{\mathfrak g}(a_i\pt a_j)\pt a_1\pt\dots\hat\imath\dots\hat\jmath\dots\pt a_n\Big).
$$

On pose $b_1=\ell_{\mathfrak g}(a_i\pt a_j),b_2=a_1, \dots, b_{n-1}=a_n$. Alors, on a
\begin{align*}
&F\circ\ell_{\mathfrak g}(a_1\otimes\dots\otimes a_n)=\\
=&\sum_{i<j}\varepsilon_{\deg(a)}\left(\begin{smallmatrix}1\ \dots\ n\\ i,j,1\dots\hat{\imath}\hat{\jmath}\dots n\end{smallmatrix}\right)\sum_{\ell>0}
\frac{1}{\ell!}\sum_{\begin{smallmatrix}I_1\sqcup\cdots\sqcup I_\ell=\{1,\ldots,n-1\}\\ I_1\dots I_\ell\neq\emptyset\end{smallmatrix}}\varepsilon_{\deg(b)}\left(\begin{smallmatrix}1\ \dots\ n-1\\ I_1\dots I_\ell\end{smallmatrix}\right)F_{|I_1|}(b_{I_1})\pt\dots\pt F_{|I_\ell|}(b_{I_\ell})\\
=&\sum_{i<j}\varepsilon_{\deg(a)}\left(\begin{smallmatrix}1\ \dots\ n\\ i,j,1\dots\hat{\imath}\hat{\jmath}\dots n\end{smallmatrix}\right)\sum_{\ell>0}
\frac{1}{\ell!}\sum_{\begin{smallmatrix}I_1\sqcup\cdots\sqcup I_\ell=\{1,\ldots,n-1\}\\ I_1\dots I_\ell\neq\emptyset\\b_1\in I_k\end{smallmatrix}}\varepsilon_{\deg(b)}\left(\begin{smallmatrix}1\ \dots\ n-1\\ I_k I_1\dots I_\ell\end{smallmatrix}\right)F_{|I_k|}(b_{I_k})\pt F_{|I_1|}(b_{I_1})\pt\dots\pt F_{|I_\ell|}(b_{I_\ell})\\
=&\sum_{i<j}\varepsilon_{\deg(a)}\left(\begin{smallmatrix}1\ \dots\ n\\ i,j,1\dots\hat{\imath}\hat{\jmath}\dots n\end{smallmatrix}\right)\sum_{\ell>0}
\frac{1}{\ell!}\sum_{\begin{smallmatrix}I_1\sqcup\cdots\sqcup I_\ell=\{1,\ldots,n-1\}\\ I_1\dots I_\ell\neq\emptyset\\b_1\in I_k, |I_k|\leq n-2\end{smallmatrix}}\varepsilon_{\deg(b)}\left(\begin{smallmatrix}1\ \dots\ n-1\\ I_k I_1\dots I_\ell\end{smallmatrix}\right)F_{|I_k|}(b_{I_k})\pt F_{|I_1|}(b_{I_1})\pt\dots\pt F_{|I_\ell|}(b_{I_\ell})\\
&\hskip 1cm+\sum_{i<j}\varepsilon_{\deg(a)}\left(\begin{smallmatrix}1\ \dots\ n\\ i,j,1\dots\hat{\imath}\hat{\jmath}\dots n\end{smallmatrix}\right)F_{n-1}(\ell_{\mathfrak g}(a_i\pt a_j)a_1\pt\dots \pt a_n).
\end{align*}
D'o\`u,
\begin{align*}
&\ell_{\mathfrak h}\circ F(a_1\otimes\dots\otimes a_n)-F\circ\ell_{\mathfrak g}(a_1\otimes\dots\otimes a_n)=0=\\
=&\sum_{j>0}\frac{1}{j!}\hskip-0.3cm\sum_{\begin{smallmatrix}I_1\sqcup\cdots\sqcup I_j=\{1,\ldots,n\}\\ I_1\dots I_j\neq\emptyset\end{smallmatrix}}\sum_{s\neq k\atop |I_s|=1, |I_k|\leq n-2}\hskip-0.5cm\varepsilon_{\deg(a)}\left(\begin{smallmatrix}1,\dots,n\\ s,I_k,I_1,\dots \hat{sk}\dots,I_j\end{smallmatrix}\right)\ell_{\mathfrak h}\Big(a_s,F_{|I_k|}(a_{I_k})\Big)F_{|I_1|}(a_{I_1})\pt\dots\pt F_{|I_j|}(a_{I_j})\\
&-\sum_{i<j}\varepsilon_{\deg(a)}\left(\begin{smallmatrix}1\ \dots\ n\\ i,j,1\dots\hat{\imath}\hat{\jmath}\dots n\end{smallmatrix}\right)\sum_{\ell>0}
\frac{1}{\ell!}\hskip-0.3cm\sum_{\begin{smallmatrix}I_1\sqcup\cdots\sqcup I_\ell=\{1,\ldots,n-1\}\\ I_1\dots I_\ell\neq\emptyset\\b_1\in I_k, |I_k|\leq n-2\end{smallmatrix}}\varepsilon_{\deg(b)}\left(\begin{smallmatrix}1\ \dots\ n-1\\ I_k I_1\dots I_\ell\end{smallmatrix}\right)F_{|I_k|}(b_{I_k})\pt F_{|I_1|}(b_{I_1})\pt\dots\pt F_{|I_\ell|}(b_{I_\ell})\\
&+\sum_{j>0}\frac{1}{j!}\hskip-0.3cm\sum_{\begin{smallmatrix}I_1\sqcup\cdots\sqcup I_j=\{1,\ldots,n\}\\ I_1\dots I_j\neq\emptyset\end{smallmatrix}}\sum_{s} \varepsilon_{\deg(a)}\left(\begin{smallmatrix}1,\dots,n\\ s,1,\dots \hat{s}\dots,n\end{smallmatrix}\right)\ell_{\mathfrak h}\Big(a_s,F_{n-1}(a_1\otimes\dots\hat{s}\dots\otimes a_{n})\Big)\\
&-\sum_{i<j}\varepsilon_{\deg(a)}\left(\begin{smallmatrix}1\ \dots\ n\\ i,j,1\dots\hat{\imath}\hat{\jmath}\dots n\end{smallmatrix}\right)F_{n-1}(\ell_{\mathfrak g}(a_i\pt a_j)a_1\pt\dots \pt a_n)\\
=&(I)-(II)+(III)-(IV).\end{align*}

Sachant que $\partial c_{|I_k|}=0$, pour tout $I_k$ tel que $|I_k|\leq n-2$, on obtient $(I)-(II)=0$. Il restera donc $(III)-(IV)=\partial c_{n-1}(a_1\otimes\dots\otimes a_n)=0$.

\

\noindent
3. Si $c=\ell_{\mathfrak h}\circ(\iota\bullet b)-b\circ\ell_{\mathfrak g} $, alors,
$$\aligned
c_{n-1}(a_0\pt\dots\pt a_n)=&\ell_{\mathfrak h}(\sum_i(-1)^{\deg(a_i)}\varepsilon_{\deg(a)}\left(\begin{smallmatrix}1~~~\dots~~~ n+1\\
i,1,\dots\hat\imath\dots n+1\end{smallmatrix}\right)
a_i\pt b(a_0\pt\dots\hat{\imath}\dots\pt a_n)\\
&-\sum_{i<j}\varepsilon_{\deg(a)}\left(\begin{smallmatrix}
0\ \dots \n\\
i,j,0\dots\hat{\imath}\hat{\jmath}\dots n
\end{smallmatrix}\right)
b_{n-2}(\ell_{\mathfrak g}(a_i\pt a_j)\pt a_0\pt\dots\hat{\imath}\dots\hat{\jmath}\dots\pt a_n)\\
=&\partial b_{n-2}(a_0\pt\dots\pt a_n).
\endaligned
$$


\section{R\'esolution de Koszul pour une alg\`ebre de Lie}\label{sec5}

On peut, comme pour l'homologie et la cohomologie de Hochschild, d\'ecrire les complexes d'homologie et de cohomologie de Chevalley comme des complexes d\'eriv\'es d'une r\'esolution, celle du module trivial $\R$. Introduisons l'alg\`ebre enveloppante ${\mathcal U}(\mathfrak g)$ de l'alg\`ebre de Lie $\mathfrak g$. On rappelle, que si $(a_i)_{i\in I}$ est une base de $\mathfrak g$, form\'ee d'\'el\'ements homog\`enes, si $I$ est totalement ordonn\'e, alors les mon\^omes $a_{i_1}a_{i_2}\dots a_{i_p}$ (calcul\'es dans $\mathcal U(\mathfrak g)$) o\`u $i_k\leq i_{k+1}$ et $i_k=i_{k+1}$ implique $a_{i_k}$ pair forment une base de $\mathcal U(\mathfrak g)$, dite base de Poincar\'e-Birkhof-Witt. On peut prouver ce r\'esultat en reprenant la preuve de \cite{[B]} donn\'ee dans le cas non gradu\'ee et en l'adaptant au cas gradu\'e.\\

Rappelons que de m\^eme une base de $S(\mathfrak g)$ est donn\'ee par les m\^emes mon\^omes, calcul\'es dans $S(\mathfrak g)$. Enfin une base de $\bigwedge\mathfrak g$ (qui est lin\'eairement isomorphe \`a $S(\mathfrak g[1])$) est donn\'ee par les mon\^omes $a_{i_1}\wedge a_{i_2}\wedge\dots\wedge a_{i_p}$, o\`u $i_k\leq i_{k+1}$ et $i_k=i_{k+1}$ implique $a_{i_k}$ impair, calcul\'es dans $\bigwedge\mathfrak g$.\\

Consid\`erons la suite exacte~:
$$
0\stackrel{\partial_{-1}}{\longleftarrow}\R\stackrel{\partial_0}{\longleftarrow}{\mathcal U}(\mathfrak g)\stackrel{\partial_1}{\longleftarrow}{\mathcal U}(\mathfrak g)\otimes\mathfrak g\stackrel{\partial_2}{\longleftarrow}\dots \stackrel{\partial_n}{\longleftarrow}{\mathcal U}(\mathfrak g)\otimes\wedge^n\mathfrak g\stackrel{\partial_{n+1}}{\longleftarrow}\dots
$$
que l'on note
$$
0\stackrel{\partial_{-1}}{\longleftarrow}\R\stackrel{\partial_0}{\longleftarrow}C_1\stackrel{\partial_1}{\longleftarrow}C_2
\stackrel{\partial_2}{\longleftarrow}\dots \stackrel{\partial_n}{\longleftarrow}C_n\stackrel{\partial_{n+1}}{\longleftarrow}\dots
$$
o\`u $\partial_0$ est l'augmentation de ${\mathcal U}(\mathfrak g)$ et $\partial_n$ est un morphisme de ${\mathcal U}(\mathfrak g)$-module d\'efini par~:
$$\aligned
\partial_n(u\otimes a_1\wedge&\dots\wedge a_n)=\sum_{i=1}^n(-1)^{i+1}\varepsilon_{|a|}(i,1\dots\hat{\imath}\dots,n)ua_i\otimes a_1\wedge\dots\hat{\imath}\dots\wedge a_n\cr
&+\sum_{i<j}(-1)^{i+j}\varepsilon_{|a|}(i,j,1,\dots\hat{\imath}\hat{\jmath}\dots,n)u\otimes[a_i,a_j]\wedge a_1\wedge\dots\hat{\imath}\dots\hat{\jmath}\dots\wedge a_n.\endaligned
$$
Cette suite exacte est une r\'esolution, elle s'appelle la r\'esolution de Koszul du $\mathfrak g$ mo\-dule trivial.\\

Pour montrer que la suite est une r\'esolution, on filtre d'abord l'espace total
$$
C=\oplus_{n\geq0}C_n=\oplus_{n\geq0}{\mathcal U}(\mathfrak g)\otimes\wedge^n\mathfrak g
$$
en posant~:
$$
F_p(C)=\bigoplus_{k=0}^p{\mathcal U}(\mathfrak g)_{p-k}\otimes\wedge^k\mathfrak g=\bigoplus_{k=0}^pF_p(C)_k.
$$

On pose $W_p^q=F_p(C)_q/F_{p-1}(C)_q$.
En ne retenant que la premi\`ere partie de l'op\'erateur $\partial_n$, on pose~:
$$
d_n(u\otimes a_1\wedge\dots\wedge a_n)=\sum_{i=1}^n(-1)^{i+1}\varepsilon_{|a|}(i,1\dots\hat{\imath}\dots,n)ua_i\otimes a_1\wedge\dots\hat{\imath}\dots\wedge a_n.
$$
Cet op\'erateur passe au quotient et d\'efinit donc une famille de suites exactes, pour $p\geq0$~:
$$
0\stackrel{d_{-1}}{\longleftarrow}\R\stackrel{d_0}{\longleftarrow}W_p^0\stackrel{d_1}{\longleftarrow}W_p^1\stackrel{d_2}{\longleftarrow}\dots \stackrel{d_p}{\longleftarrow}W_p^p\stackrel{d_{p+1}}{\longleftarrow}0.
$$
Ces suites admettent une homotopie qui est d\'efinie sur la base de $F_k(C)_n$~:
$$
h_{n+1}(a_{i_1}\dots a_{i_k}\otimes a_{j_1}\wedge\dots\wedge a_{j_n})=\left\{\aligned&0~~\text{ si }~~i_k<j_n~~\text{ ou si }~~i_k=j_n~~\text{ et }~~a_{i_k}~~\text{ est pair}\\
&\frac{1}{\#\{\ell,j_\ell=j_n\}}a_{i_1}\dots a_{i_{k-1}}\otimes a_{i_k}\wedge a_{j_1}\wedge\dots\wedge a_{j_n}~~\text{ sinon}.\endaligned\right.
$$
On v\'erifie que $h_n$ passe au quotient, il est d\'efini sur $W_k^n$ et pour tout $k$ et tout $n$,
$$
h_n\circ d_n+d_{n+1}\circ h_{n+1}=id_{W_k^n}
$$

Maintenant on a un diagramme commutatif pour tout $p$~:

\noindent
$
0\longleftarrow \mathbb{R}\longleftarrow \hskip 0.5cm W_p^0 \hskip 0.5cm \longleftarrow
\hskip 0.5cm W_p^1\hskip 0.5cm \longleftarrow\cdots\longleftarrow\hskip 0.5cm
W_p^q\hskip 0.5cm \longleftarrow\cdots\longleftarrow \hskip 0.3cm W_p^p\hskip 0.4cm \longleftarrow0
$\vskip 0.2cm

\hskip 0.7cm$\parallel\hskip 1.1cm\pi_0\uparrow\hskip
1.6cm\pi_1\uparrow\hskip 3.2cm\pi_q\uparrow\hskip
2.9cm\pi_p\uparrow$\vskip 0.2cm

\noindent$ 0\longleftarrow \mathbb{R}\longleftarrow \hskip 0.2cm
(F_pC)_0\hskip 0.1cm  \longleftarrow \hskip 0.2cm (F_pC)_1\hskip
0.2cm \longleftarrow\cdots\longleftarrow\hskip 0.2cm (F_pC)_q\hskip
0.2cm \longleftarrow\cdots\longleftarrow (F_pC)_p\longleftarrow0
$\vskip 0.2cm

\hskip 0.7cm$\parallel\hskip 1.2cm j_0\uparrow\hskip 1.7cm
j_1\uparrow\hskip 3.2cm j_q\uparrow\hskip 2.9cm j_p\uparrow$\vskip
0.2cm

\noindent$ 0\longleftarrow \mathbb{R}\longleftarrow
(F_{p-1}C)_0\longleftarrow
(F_{p-1}C)_1\longleftarrow\cdots\longleftarrow
(F_{p-1}C)_q\longleftarrow\cdots\longleftarrow (F_{p-1}C)_p=0 $

\vskip 0.4cm
\noindent
ici, $j_q$ est l'injection canonique et $\pi_q$ est la projection canonique et les fl\`eches de la premi\`ere ligne sont les $d_k$, celles des lignes 2 et 3 sont les $\partial_k$.

On v\'{e}rifie directement que ces diagrammes commutent et on vient de voir que la premi\`ere ligne est exacte.

On en d\'{e}duit de fa\c con classique une suite longue en homologie $\forall p\geq1$:

$$
\cdots\longrightarrow H((F_{p-1}C)_q)\longrightarrow H((F_{p}C)_q)\longrightarrow H(W^q_p)=0\longrightarrow
H((F_{p-1}C)_{q+1})\longrightarrow\cdots
$$
Alors, pour tout $p\geq1$ et tout $q$, $H((F_{p-1}C)_q)=H((F_{p}C)_q)$.

De plus, on a $(F_0C)_0=\mathbb{R}$ et $(F_0C)_{-1}=\mathbb{R}$. Donc, on obtient $H((F_0C)_q)=0$ pour tout $q$, puis on en d\'{e}duit que
$$
H((F_{p}C)_q)=0
$$
pour tout $p$ et tout $q$. Ou\\

\begin{prop}{\rm(La r\'esolution de Koszul est une r\'esolution)}

\

Le complexe de Koszul du $\mathfrak g$ module trivial $\R$ est une r\'esolution.\\
\end{prop}

Soit $M$ un $\mathfrak g$ module \`a gauche, on le consid\`ere comme un $\mathcal U(\mathfrak g)$ module \`a gauche. Alors, apr\`es d\'ecalage, on a les identifications d'espaces vectoriels~:
$$
M\otimes_{\mathcal U(\mathfrak g)}C_n\simeq M[1]\otimes S^n(\mathfrak g[1])\quad\text{ et }\quad Hom_{\mathcal U(\mathfrak g)}(C_n,M)\simeq L(S^n(\mathfrak g[1]),M[1]).
$$

L'homologie et la cohomologie de Chevalley du module $M$, d\'efinies dans la section pr\'ec\'edente sont alors obtenues comme des complexes d\'eriv\'es de la r\'esolution de Koszul par un argument du m\^{e}me type que celui de la section 1.\\

\section{$G_\infty$ alg\`ebres}\label{sec6}

Consid\'erons maintenant une alg\`ebre de Gerstenhaber $G$. $G$ est un espace vectoriel gradu\'e, muni d'une multiplication associative et commutative, de degr\'e 0, not\'ee $(a,b)\mapsto ab$ et d'un crochet not\'e $(a,b)\mapsto[a,b]$, de degr\'e -1 tel que $(G[1],[~,~])$ est une alg\`ebre de Lie gradu\'ee. De plus l'application $ad_a~:b\mapsto[a,b]$ est une d\'erivation de degr\'e $|a|-1$ pour la multiplication (relation de Leibniz).\\

La d\'efinition de l'alg\`ebre $G_\infty$ enveloppante comme ci-dessus, demande la construction d'une cod\'erivation, combinaison du produit et du crochet, telle que la relation de structure, c'est \`a dire l'associativit\'e du produit, la relation de Jacobi pour le crochet et la relation de Leibniz soient \'equivalentes \`a l'annulation du carr\'e  de cette cod\'erivation. La premi\`ere difficult\'e est le fait que le produit et le crochet ne sont pas de m\^eme degr\'e et qu'ils d\'ependent du m\^eme nombre d'arguments. Si on remplace $G$ par $G[k]$, le produit et le crochet restent de degr\'es diff\'erents. La seconde difficult\'e vient des relations de sym\'etries diff\'erentes pour chaque loi.\\

La solution propos\'ee par Ginot \cite{[G]} consiste \`a remplacer le produit par la cod\'erivation $m$ comme pour les alg\`ebres commutatives, \`a prolonger le crochet en un crochet sur la $C_\infty$ alg\`ebre enveloppante de $G$. D\`es lors cette alg\`ebre poss\`edera deux op\'erations, l'une ayant un seul argument ($m$) de degr\'e 1 l'autre \'etant un crochet de degr\'e 0 entre deux arguments. Un nouveau d\'ecalalge et la construction d'une $L_\infty$ alg\`ebre enveloppante de $C_\infty(G)$ permet de se retrouver avec deux cod\'erivations de degr\'e 1 que l'on peut additionner et ainsi de d\'ecrire la structure comme l'annulation d'un carr\'e.\\

Si on ne regarde que les sym\'etries des relations, on trouve~:
$$\aligned
&Com(a_1\otimes a_2)=a_1a_2-(-1)^{|a_1||a_2|}a_2a_1=0,\\
&Ass(a_1\otimes a_2\otimes a_3)=a_1(a_2a_3)-(a_1a_2)a_3=0,\\
&Antisym(a_1\otimes a_2)=[a_1,a_2]+(-1)^{\deg(a_1)\deg(a_2)}[a_2,a_1],\\
&Jac(a_1\otimes a_2\otimes a_3)=\oint_{1,2,3}(-1)^{\deg(a_1)\deg(a_3)}\big[[a_1,a_2],a_3\big]=0\\
&Leibn(a_1\otimes a_2\otimes a_3)=[a_1,a_2a_3]-[a_1,a_2]a_3-(-1)^{|a_2|\deg(a_1)}a_2[a_1,a_3].
\endaligned
$$
($\oint_{123}$ signifie somme sur les permutations circulaires sur 1,2,3). De ces relations, on peut construire une relation qui ne v\'erifie pas la r\`egle des signes de Koszul~:
$$
[a_1a_2,a_3]=a_1[a_2,a_3]+(-1)^{|a_2|\deg(a_3)}[a_1,a_3]a_2.
$$
On va donc traiter en deux \'etapes la construction de la cog\`ebre qu'on utilisera. D'abord, on construit la cog\`ebre de Lie
$$
(\mathcal H,\delta)=(\oplus_{n\geq1}\underline{\otimes^nG[1]},\delta),
$$
cette cog\`ebre est d\'efinie dans la section \ref{sec4}.

Ensuite on \'etend le cocrochet $\delta$ en un cocrochet $\kappa$ de $\mathcal H[1]$ en posant,
$$
X=\underline{a_1\otimes\dots\otimes a_p},
$$
et
$$\aligned
\kappa(X)&=\sum_{j=1}^{p-1}\Big((-1)^{\sum_{k\leq j}\deg(a_k)}\underline{a_1\otimes\dots\otimes a_j}\bigotimes\underline{
a_{j+1}\otimes a_p}\cr
&-\varepsilon_{\deg(a)}\left(\begin{smallmatrix}1,\dots,\dots,p\\
j+1,\dots,p,1,\dots,j\end{smallmatrix}\right)(-1)^{\sum_{k>j}\deg(a_k)}\underline{a_{j+1}\otimes\dots\otimes
a_p} \bigotimes\underline{a_1\otimes\dots\otimes
a_j}\Big).\endaligned
$$

Cette formule s'\'ecrit de fa\c con condens\'ee ainsi~: posons 
$$
U_j=\underline{a_1\otimes\dots\otimes a_j},\quad V_j=\underline{a_{j+1}\otimes\dots\otimes a_p},
$$
ce sont des \'el\'ements de $\mathcal H[1]$ de degr\'es 
$$
\deg(U_j)=\sum_{k\leq j}\deg(a_k)-1~~\text{ et }~~\deg(V_j)=\sum_{k>j}\deg(a_k)-1.
$$
Alors~:
$$
\kappa(X)=\sum_{j=1}^{p-1}(-1)^{\deg(U_j)+1}\Big(U_j\bigotimes V_j+\tau\big(U_j\bigotimes V_j\big)\Big).
$$

Le cocrochet $\kappa$ sur $(\underline{\otimes^pG[1]})[1]$ est alors cosym\'etrique ($\kappa=\tau\circ\kappa$ o\'u $\tau$ est la volte dans $\mathcal H[1]$) et de degr\'e 1.\\

Comme dans la section \ref{sec5}, on consid\`ere la cog\`ebre $(S^+(\mathcal H[1]),\Delta)=(\oplus_{n\geq1}S^n(\mathcal H[1]),\Delta)$, qui traduira les sym\'etries des relations $Jac$ et $Antisym$.

On prolonge alors $\kappa$ \`a $S^+(\mathcal{H}[1])$ en posant~:
$$\aligned
\kappa(X_1\pt\dots\pt X_n)=&\sum_{\begin{smallmatrix}1\leq s\leq n\\ I\cup J=\{1,\dots,n\}\setminus\{s\}\end{smallmatrix}}(-1)^{\sum_{i<s}x_i}\sum_{\begin{smallmatrix}U_s\otimes V_s=X_s\\ U_s,V_s\neq\emptyset\end{smallmatrix}} (-1)^{\deg(U_s)+1}\times\cr
&\times\left(\varepsilon\left(\begin{smallmatrix} x_1,\dots,x_n\\ x_I~u_s~v_s~x_J\end{smallmatrix}\right)X_I\pt U_s\bigotimes V_s \pt X_J+\varepsilon\left(\begin{smallmatrix} x_1\dots x_n\\ x_I~v_s~u_s~x_J\end{smallmatrix}\right)X_I\pt V_s\bigotimes U_s \pt X_J\right).
\endaligned
$$

Les notations sont celles de \cite{[AAC3]}~:
$$
x_j=\deg(X_j),~~X_j\in\mathcal H[1],~~x_I=(deg(x_k),~k\in I),~~u_s=\deg(U_s),~~v_s=\deg(V_s)
$$ 
et
$\varepsilon\left(\begin{smallmatrix} x_1,\dots,x_n\\ x_I~u_s~v_s~x_J\end{smallmatrix}\right)$ repr\'esente le signe de la permutation $(1,\dots,n)\mapsto I\cup\{s\}\cup J$, en tenant compte des degr\'es, $\varepsilon\left(\begin{smallmatrix} x_1,\dots,x_n\\ x_I~v_s~u_s~x_J\end{smallmatrix}\right)$ est 
$$
(-1)^{\deg(U_s)\deg(V_s)}\times\text{ ce signe}.
$$
On rappelle que $\deg(X_s)=\deg(U_s)+\deg(V_s)+1$.

Rappelons que $\Delta$ est cocommutative et coassociative, de degr\'e 0~:\\

(i) $\tau\circ\Delta=\Delta$ (cocommutativit\'e),\\

(ii) $(id\otimes\Delta)\circ\Delta=(\Delta\otimes id)\circ\Delta$ (coassociativit\'e).\\

En posant $\tau_{12}=\tau\otimes id$ et $\tau_{23}=id\otimes\tau$, $\kappa$ est de degr\'e 1, cocommutatif et v\'erifie les identit\'es de coJacobi et de coLeibniz suivantes~:\\

(iii) $\tau\circ\kappa=\kappa$ (cocommutativit\'e),\\

(iv) $\Big(id^{\otimes3}+\tau_{12}\circ\tau_{23}+\tau_{23}\circ\tau_{12}\Big)\circ(\kappa\otimes id)\circ\kappa=0$ (identit\'e de coJacobi),\\

(v) $(id\otimes\Delta)\circ\kappa=(\kappa\otimes id)\circ\Delta+\tau_{12}\circ(id\otimes\kappa)\circ\Delta$ (identit\'e de coLeibniz).\\

(Voir \cite{[BGHHW]} et \cite {[AAC3]}).\\

L'espace $S^+(\mathcal{H}[1])$ est maintenant une bicog\`ebre $(S^+(\mathcal{H}[1]),\kappa,\Delta)$. Cette structure est libre dans le sens que tout morphisme $F$ ou toute cod\'erivation $D$ peut \^etre reconstruite \`a partir de leur s\'erie de Taylor~:
$$\aligned
F_{p_1,\dots,p_r}~&:(\underline{\otimes^{p_1}G[1]})[1]\pt(\underline{\otimes^{p_2}G[1]})[1]\pt\dots\pt(\underline{\otimes^{p_r}G[1]})[1] \longrightarrow (G'[1])[1]\\
D_{p_1,\dots,p_r}~&:(\underline{\otimes^{p_1}G[1]})[1]\pt(\underline{\otimes^{p_2}G[1]})[1]\pt\dots\pt(\underline{\otimes^{p_r}G[1]})[1] \longrightarrow (G[1])[1]
\endaligned
$$

Comme on a vu comment construire des homomorphismes ou des cod\'erivations de $S^+\mathcal H[1]$ \`a partir de leur s\'erie de Taylor, pour reconstruire $F$ et $D$ \`a partir des $F_{p_1,\dots,p_r}$ ou des $D_{p_1,\dots,p_r}$, il nous suffit de d\'efinir les applications~:
$$
F_n:\oplus_{n\geq1}S^n(\mathcal H[1])\longrightarrow\mathcal H'[1]\quad\text{ et }\quad
D_n:\oplus_{n\geq1}S^n(\mathcal H[1])\longrightarrow\mathcal H[1].
$$

Soit donc $X_1\pt\dots\pt X_n$ un \'el\'ement de $(\underline{\otimes^{p_1}G[1]})[1]\pt(\underline{\otimes^{p_2}G[1]})[1]\pt\dots \pt (\underline{\otimes^{p_r}G[1]})[1]$, avec
$$
X_j=\underline{a_1^j\otimes\dots\otimes a^j_{p_j}}.
$$
$F_n(X_1\pt\dots\pt X_n)$ (resp. $D_n(X_1\pt\dots\pt X_n)$) est une somme de produits tensoriels modulo les battements de $F_{q_1,\dots,q_s}(Y_k)$ ($1\leq k\leq t$) (resp. de $D_{q_1,\dots,q_s}(Y_k)$ et de $Y_k$) o\`u les $Y_k$ sont des produit $\pt$ de parties $U_i^j$ des $X_j$ de la forme~:
$$
U_i^j=\underline{a^j_{r^j_i+1}\otimes a^j_{r^j_i+2}\otimes\dots\otimes a^j_{r^j_{i+1}}}\qquad(0\leq r^j_s\leq p_j).
$$

On consid\`ere toutes les d\'ecompositions possibles des $X_j$ en produit $\underline\otimes$ de $U_i^j$, on permute les $U_i^j$, on note $V_k^\ell$ les $U_i^j$ apr\`es cette permutation, on pose $Y_k=V_1^k\pt\dots\pt V_{s_k}^k$. On note
$$
V^k_\ell\in(\underline{\otimes^{q^k_\ell}G[1]})[1].
$$

Les $V^k_\ell$ forment une permutation des $U_i^j$. Si un $X_j$ n'est pas d\'ecompos\'e ($r^j_j=1$), il ne peut appara\^\i tre qu'en facteur d'au moins une vraie partie $U^{j'}_i$ d'un autre $X$ ($r^{j'}_{j'}>1$). Si un $X_j$ est d\'ecompos\'e ($r^j_j>1$) chacune de ses parties appara\^\i t dans un $Y$ diff\'erent enfin il y a autant de $\pt$ et de $\underline\otimes$ dans l'expression
$$
X_1\pt\dots\pt X_n=\left(\underline{U^1_1\otimes\dots\otimes U^1_{r^1_1}}\right)\pt\dots\pt\left(\underline{U^n_1\otimes\dots\otimes U^n_{r^n_n}}\right)
$$
que dans l'expression (formelle)~:
$$
\underline{Y_1\otimes\dots\otimes Y_t}=\underline{(V_1^1\pt\dots\pt V^1_{s_1})\otimes\dots\otimes(V^t_1\pt\dots\pt V_{s_t}^t)}.
$$
(Voir \cite{[AAC3]}).

$F_n$ envoie le produit $X_1\pt\dots\pt X_n$ sur des sommes de termes de la forme
$$\aligned
F_n(X_1\pt\dots\pt X_n)&=\sum_{U,V}\pm \underline{F_{q^1_1,\dots,q^1_{s_1}}(Y_1)\otimes\dots\otimes F_{q^t_1,\dots,q^t_{s_t}}(Y_t)}\\
&=\sum_{U,V}\pm \underline{F_{q^1_1,\dots,q^1_{s_1}}(V_1^1\pt\dots\pt V^1_{s_1})\otimes\dots\otimes F_{q^t_1,\dots,q^t_{s_t}}(V^t_1\pt\dots\pt V_{s_t}^t)}.
\endaligned
$$
De m\^eme $D_n$ envoie $X_1\pt\dots\pt X_n$ sur des sommes de termes de la forme
$$\aligned
D_n(X_1\pt\dots\pt X_n)&=\sum_{U,V}\pm \underline{Y_1\otimes\dots\otimes D_{q^k_1,\dots,q^k_{s_k}}(Y_k)\otimes\dots\otimes Y_t}\\
&=\sum_{U,V}\pm \underline{Y_1\otimes\dots\otimes D_{q^k_1,\dots,q^k_{s_k}}(V^k_1\pt\dots\pt V_{s_k}^k)\otimes\dots\otimes Y_t}.
\endaligned
$$

On va expliciter cette construction pour les lois de $G$, que l'on \'etend \`a la bicog\`ebre $(S^+(\mathcal{H}[1]),\kappa,\Delta)$.\\

\begin{prop}{\rm(Construction de $m$ et $\ell$)}

\

Soit $(G,.,[~,~])$ une alg\`ebre de Gerstenhaber. Le produit commutatif s'\'etend \`a $\oplus_{n\geq1}\underline{\otimes^n G[1]}$ comme ci-dessus par~:
$$\aligned
m(\underline{a_1\otimes a_2})&=(-1)^{\deg(a_1)}a_1a_2\\
m(\underline{a_0\otimes\dots\otimes a_n})&=\sum_{j=0}^{n-1}(-1)^{\sum_{i<j}\deg(a_i)}\underline{a_1\otimes\dots\otimes m(a_j\otimes a_{j+1})\otimes\dots\otimes a_n}\quad(a_i^j\in G[1])\\
m(X_1\pt\dots\pt X_n)&=\sum_{j=1}^n(-1)^{\sum_{k<j}\deg(X_k)}X_1\pt\dots\pt m(X_j)\pt\dots\pt X_n\quad (X_k\in\mathcal H[1]).
\endaligned
$$

Ainsi d\'efini, $m$ est une cod\'erivation de $\kappa$ et de $\Delta$.

Le crochet de Lie s'\'etend d'abord en un crochet sur $\mathcal H$, ainsi~:
$$\aligned
&\left[\underline{a_1\otimes\dots\otimes a_p}\wedge\underline{a_{p+1}\otimes\dots\otimes a_{p+q}}\right]=\\
&\hskip 1cm\sum_{\begin{smallmatrix}\sigma\in Bat(p,q)\\ \sigma^{-1}(k)\leq p<\sigma^{-1}(k+1)\end{smallmatrix}}
\hskip-0.5cm \varepsilon_{\deg(a)}(\sigma^{-1})\underline{a_{\sigma^{-1}(1)}\otimes\dots\otimes[a_{\sigma^{-1}(k)},a_{\sigma^{-1}(k+1)}]
\otimes\dots\otimes a_{\sigma^{-1}(p+q)}}.
\endaligned
$$
Enfin on \'etend ce crochet \`a $S^+(\mathcal H[1])$ comme ci-dessus par~:
$$\aligned
\ell(X_1\pt X_2)&=(-1)^{\deg(X_1)}\left[X_1\wedge X_2\right]\\
\ell(X_1\pt\dots\pt X_n)&=\sum_{i<j}\varepsilon_{\deg(X)}\left(\begin{smallmatrix}1~\dots ~\dots~ n\\ i,j,1,\dots\hat{\imath}\hat{\jmath}\dots,n\end{smallmatrix}\right) \ell(X_i\pt X_j)\pt X_1\dots\hat{\imath}\dots\hat{\jmath}\dots X_n.
\endaligned
$$

Ainsi d\'efini, $\ell$ est une cod\'erivation de $\kappa$ et de $\Delta$.

La strucure d'alg\`ebre de Gerstenhaber, c'est \`a dire l'associativit\'e du produit, la relation de Jacobi pour le crochet et la relation de compatibilit\'e de Leibniz, est \'equivalente \`a
$$
(m+\ell)\circ(m+\ell)=0.
$$
\end{prop}

On a donc les propri\'et\'es suivantes~: $m$ est de degr\'e 1 et~:\\

(i) $(m\otimes id+id\otimes m)\circ\kappa=-\kappa\circ m$ ($\kappa$-cod\'erivation gradu\'ee)\\

(ii) $(m\otimes id+id\otimes m)\circ\Delta=\Delta\circ m$ ($\Delta$-cod\'erivation)\\

(iii) $m\circ m=0$ (associativit\'e)\\

$\ell$ est de degr\'e 1 et~:\\

(iv) $(\ell\otimes id+id\otimes\ell)\circ\kappa=-\kappa\circ\ell$ ($\kappa$-cod\'erivation gradu\'ee)\\

(v) $(\ell\otimes id+id\otimes\ell)\circ\Delta=\Delta\circ\ell$ ($\Delta$-cod\'erivation)\\

(vi) $\ell\circ\ell=0$ (relation de Jacobi)\\

Enfin~:\\

(vii) $\ell\circ m+m\circ\ell=0$ (relation de Leibniz).\\

On pose donc~:\\

\begin{defn}{\rm(Alg\`ebre $G_\infty$ enveloppante)}

\

Soit $(G,.,[~,~])$ une alg\`ebre de Gerstenhaber gradu\'ee. On appelle $G_\infty$ alg\`ebre enveloppante de $G$ la bicog\`ebre de Lie , coassociative et cocommutative
$$
(\oplus_{n\geq1}S^n(\oplus_{p\geq1}(\underline{\otimes^pG[1]})[1]),\kappa,\Delta).
$$
munie de la cod\'erivation $m+\ell$. On note cette alg\`ebre enveloppante
$$
G_\infty(G)=(\oplus_{n\geq1}S^n(\oplus_{p\geq1}(\underline{\otimes^p G[1]})[1]),\kappa,\Delta,m+\ell).
$$
\end{defn}

La $G_\infty$ alg\`ebre enveloppante de $G$ est donc un complexe d'homologie. Posons
$$
C_N(G)=\sum_{r=1}^N\sum_{p_1+\dots+p_r=N}\left(\underline{\otimes^{p_1}G[1]}\right)[1]\pt\dots\pt\left(\underline{\otimes^{p_r}G[1]}\right)[1],
$$
le complexe est~:
$$
0\stackrel{m+\ell}{\longleftarrow}C_1(G)\stackrel{m+\ell}{\longleftarrow}C_2(G)\stackrel{m+\ell}{\longleftarrow}\dots \stackrel{m+\ell}{\longleftarrow}C_{N-1}(G)\stackrel{m+\ell}{\longleftarrow}\dots
$$

\begin{defn}{\rm($G$-module)}

\

Soit $G$ une alg\`ebre de Gerstenhaber. Un $G$-module $M$ est un espace vectoriel gradu\'e muni de deux lois externes~:
$$
.~:G\times M\longrightarrow M,\qquad (a,m)\mapsto am,
$$
telle que $M$ devienne un bimodule sym\'etrique pour la multiplication commutative de $G$ ($ma=(-1)^{|a||m|}am$) et
$$
\ppt~:G[1]\times M[1]\longrightarrow M[1],\qquad (a,m)\mapsto a\ppt m,
$$
telle que $M$ soit un module pour la structure d'alg\`ebre de Lie de $G$. De plus les lois $.$ et $\ppt$ sont compatibles entre elles, c'est \`a dire satisfont, pour tout $a_1$, $a_2$ de $G$ et $m$ de $M$,
$$
[a_1,a_2]m=a_1\ppt(a_2m)-(-1)^{|a_2|\deg(a_1)}a_2(a_1\ppt m)
$$
et
$$
(a_1a_2)\ppt m=a_1(a_2\ppt m)+(-1)^{|m||a_2|}a_2(a_1\ppt m).
$$
\end{defn}

Soit $G$ une alg\`ebre de Gerstenhaber et $M$ un $G$-module. Imitant les constructions des sections pr\'ec\'edentes, on peut d\'efinir une nouvelle alg\`ebre de Gerstenhaber, not\'ee $G\ltimes M$ en posant $G\ltimes M=G\oplus M$ en tant qu'espace vectoriel et~:
$$
(a+u).(b+v)=ab+av+ub\qquad\text{ et }\quad[a+u,b+v]=[a,b]+a\ppt v-(-1)^{\deg(b)\deg(u)}b\ppt u.
$$
$M$ est un $G$-module si et seulement si $G\ltimes M$ est une alg\`ebre de Gerstenhaber.\\

Pour d\'efinir l'homologie du $G$-module $M$, on construit la $G_\infty$ alg\`ebre enveloppante de $G\ltimes M$
$$
0\stackrel{m+\ell}{\longleftarrow}C_1(G\ltimes M)\stackrel{m+\ell}{\longleftarrow}C_2(G\ltimes M)\stackrel{m+\ell}{\longleftarrow}\dots \stackrel{m+\ell}{\longleftarrow}C_{N-1}(G\ltimes M)\stackrel{m+\ell}{\longleftarrow}\dots
$$
et on se restreint aux sous-espaces~:
$$\aligned
C_N(G,M)&=\sum_{r=1}^N\sum_{p_1+\dots+p_r+p=N}\sum_{k=0}^p\\
&\left(\underline{\otimes^{p_1}G[1]}\right)[1]\pt\dots\pt\left(\underline{\otimes^{p_r}G[1]}\right)[1]\pt\left(\underline{\otimes^kG[1]\otimes M[1]\otimes\otimes^{p-k}G[1]}\right)[1],
\endaligned
$$
de $C_N(G\ltimes M)$. On obtient un sous-complexe qu'on appelle le complexe de l'homologie de Chevalley-Harrison du $G$-module $M$~:
$$
0\stackrel{\partial_0}{\longleftarrow}C_0(G,M)\stackrel{\partial_1}{\longleftarrow}C_1(G,M)\stackrel{\partial_2}{\longleftarrow}\dots \stackrel{\partial_{N-1}}{\longleftarrow}C_{N-1}(G,M)\stackrel{\partial_N}{\longleftarrow}\dots
$$

\begin{prop}{\rm(Expression du bord de Chevalley-Harrison)}

\

Le bord de l'homologie de Chevalley-Harrison du $G$-module $M$ est donn\'e explicitement ainsi.

Soit $Y=\underline{a_1\otimes\dots\otimes a_k\otimes v\otimes\dots\otimes a_p}$ ($a_1,\dots,a_p\in G$, $v\in M$). On rappelle que le cobord $\partial_{Har}$ de Harrison, d\'efini dans la section \ref{sec3} est $\partial_{Har}(Y)=m(Y)$.

Soit $X_i=\underline{a^i_1\otimes\dots\otimes a^i_{p_i}}$ ($a^i_1,\dots,a^i_{p_i}\in G$), on pose $X_i=\underline{b_1\otimes\dots\otimes b_{p_i}}$, $Y=\underline{b_{p_i+1}\otimes\dots\otimes b_{p_i+p+1}}$ et~:
$$
\ell(X_i,Y)=\hskip-0.7cm\sum_{\begin{smallmatrix}\sigma\in Bat(p_i,p+1)\\ \sigma^{-1}(k)\leq p_i<\sigma^{-1}(k+1)\end{smallmatrix}}
\hskip-0.5cm\varepsilon_{\deg(b)}(\sigma^{-1})b_{\sigma^{-1}(1)}\underline\otimes\dots\underline\otimes[b_{\sigma^{-1}(k)},b_{\sigma^{-1}(k+1)}]\underline\otimes\dots\underline\otimes b_{\sigma^{-1}(p+q)}.
$$
Enfin, pour $X_1\pt\dots\pt X_r\pt Y\in C_N(G,M)$,
$$\aligned
\partial_N(X_1\pt\dots\pt X_r\pt Y)&=(\ell+m)(X_1\pt\dots\pt X_r)\pt Y+(-1)^{\sum_i\deg(X_i)}(X_1\pt\dots\pt X_r)\pt m(Y)\\
&+\sum_i(-1)^{\sum_{j\neq i}\deg(X_j)} \varepsilon_{\deg(X)}\left(\begin{smallmatrix}1,\dots,\dots,r\\1,\dots\hat{\imath}\dots,r,i \end{smallmatrix}\right) X_1\pt\dots\hat{\imath}\dots\pt X_r\pt\ell(X_i,Y).
\endaligned
$$
\end{prop}

De m\^eme, la cohomologie de Chevalley-Harrison du $G$-module $M$ est d\'efinie en \'ecrivant l'\'equation de morphisme de $G_\infty$ alg\`ebre. Plus pr\'ecis\'ement, on se donne une suite d'applications $F_N:C_N(G)\longrightarrow((G\ltimes M)[1])[1]$ de la forme~:
$$
F_1=\iota+c_1,\quad F_N=c_N,
$$
o\`u $\iota$ est l'injection canonique de $G$ dans $G\ltimes M$ et les $c_N$ sont des applications lin\'eaires de $C_N(G)$ dans $M$, plus pr\'ecis\'ement,
$$
c_N=\sum_{r=1}^N\sum_{p_1+\dots+p_r=N}c_{p_1\dots p_r}
$$
avec~:
$$
c_{p_1\dots p_r}~:~\left(\underline{\otimes^{p_1}G[1]}\right)[1]\pt\dots\pt\left(\underline{\otimes^{p_r}G[1]}\right)[1]\longrightarrow M[N+r].
$$
On construit le morphisme $F$ de bicog\`ebres dont la s\'erie de Taylor est $(F_N)$~:
$$
F:(G_\infty(G),\kappa,\Delta)\longrightarrow(G_\infty(G\ltimes M),\kappa,\Delta),
$$
enfin on \'ecrit que $F$ est un morphisme pour la structure $(m_G+\ell_G)$, respectivement la structure $(m_{G\ltimes M},\ell_{G\ltimes M})$. On trouve des conditions sur les $c_N$ que l'on \'ecrit $\partial_Nc_N=0$. On obtient (voir \cite{[AAC3]})~:
$$\aligned
C^N(G,M)&=\sum_{r=1}^N\sum_{p_1+\dots+p_r=N} C^N_{p_1,\dots,p_r}\\
&=\sum_{r=1}^N\sum_{p_1+\dots+p_r=N} L\left(\left(\underline{\otimes^{p_1}G[1]}\right)[1]\pt\dots\pt \left(\underline{\otimes^{p_r}G[1]}\right)[1],M[N+r]\right),
\endaligned
$$
On obtient alors le complexe~:
$$
0\stackrel{\partial^1}{\longrightarrow}C^1(G,M)\stackrel{\partial^2}{\longrightarrow}\dots\stackrel{\partial^{N-1}}{\longrightarrow}C^{N-1}(G,M)\stackrel{\partial^N}{\longrightarrow}\dots
$$

\begin{prop}{\rm(Expression du cobord de Chevalley-Harrison de $M$)}

\

Le cobord $\partial^N:C^N(G,M)\longrightarrow C^{N+1}(G,M)$ est de la forme $\partial^N=d_m+d_\ell$ avec $d_m~:C^N_{p_1\dots p_n}\longrightarrow\displaystyle\sum_j C^{N+1}_{p_1\dots(p_j+1)\dots p_n}$, plus pr\'ecis\'ement, si
$$
X_1\pt\dots\pt X_n\in(\underline{\otimes^{p_1}G[1]})[1]\pt\dots\pt(\underline{\otimes^{p_j+1}G[1]})[1]\pt\dots\pt(\underline{\otimes ^{p_n}G[1]})[1],
$$
alors si $c^N_{p_1\dots p_n} \in C^N_{p_1\dots p_n}$, $X_j \in \otimes ^{p_j+1}G[1])[1]$, $X_i \in \otimes^{p_i}G[1])[1]$

$$
\begin{aligned}
(d_mc^N_{p_1\dots p_n})&(X_1\pt\dots\pt X_n)=\\
=&(d_mc)^{(N+1)}_{p_1\dots(p_j+1)\dots p_n}(X_1\pt\dots\pt X_n)\\
=&(-1)^{\deg(a_j^1)\sum_{i<j}\deg(X_i)} a^j_1\underline\otimes c_{p_1,\dots,p_j,\dots,p_n}^{(N)}(X_1 \pt\dots\pt(\underline{a^j_2\otimes\dots\otimes\a_{p_j+1}^j})\pt\dots\pt X_n)\\
&+(-1)^{\deg(a_j^{p_j+1})\sum_{i>j}\deg(X_i)} c_{p_1,\dots,p_j,\dots,p_r}^{(N)}(X_1\pt\dots\pt(\underline{a^j_1\otimes\dots\otimes a_{p_j}^j})\pt\dots\pt X_n)a^j_{p_j+1}\\
&-(-1)^{\sum_{i<j}\deg(X_i)}c^{(N)}_{p_1\dots p_n}(X_1\pt\dots\pt m(X_j)\pt\dots\pt X_n).
\end{aligned}
$$

De m\^eme, $d_\ell:C^N_{p_1\dots p_n}\longrightarrow\displaystyle\sum_{j,~q_1+q_2=p_j+1}C^{N+1}_{q_1,q_2,p_1\dots\hat\jmath\dots p_n}$ s'\'ecrit
$$
(d_\ell c^N_{p_1\dots p_n})=\sum_{\begin{smallmatrix}j\\ q_1+q_2=p_j+1\end{smallmatrix}}(d_\ell c)^N_{q_1,q_2,p_1\dots \hat\jmath\dots p_n}.
$$
Avec (les notations sont celles utilis\'ees dans la d\'efinition de $\kappa$)
\begin{itemize}
\item[1.] Si $q_1>1$ et $q_2>1$, alors
$$
\begin{aligned}
(d_\ell c)^{N+1}_{q_1,q_2,p_1\dots\hat\jmath\dots p_n}&(Y_1\pt Y_2\pt X_1\pt\dots\hat\jmath\dots\pt X_n)=\\
=&-\varepsilon\left(\begin{smallmatrix}y_1y_2x_1\dots\hat\jmath\dots x_n\\ x_1\dots y_1y_2\dots x_n\end{smallmatrix}\right)(-1)^{\sum_{i<j}\deg(X_i)}c^{(N)}_{p_1\dots p_n}(X_1\pt\dots\pt\ell(Y_1,Y_2)\pt\dots\pt X_n).
\end{aligned}
$$
\item[2.] Si $q_1=1$ et $q_2=p_j>1$, alors
$$
\begin{aligned}
(d_\ell c)^{N+1}_{q_1,q_2,p_1\dots\hat\jmath\dots p_n}&(Y_1\pt Y_2\pt X_1\pt\dots\hat\jmath\dots\pt X_n)=\\
=&-\varepsilon\left(\begin{smallmatrix}y_1y_2x_1\dots\hat\jmath\dots x_n\\ x_1\dots y_1y_2\dots x_n\end{smallmatrix}\right)(-1)^{\sum_{i<j}\deg(X_i)} c^{(N)}_{p_1\dots p_n}(X_1\pt\dots\pt\ell(Y_1,Y_2)\pt\dots\pt X_n)\\
&+\varepsilon\left(\begin{smallmatrix}y_1y_2x_1\dots\hat\jmath\dots x_n\\ y_1x_1\dots y_2\dots x_n\end{smallmatrix}\right) Y_1\ppt c^{(N)}_{p_1\dots p_N}(X_1\pt\dots\pt Y_2\pt\dots\pt X_n)\Big).
\end{aligned}
$$
\item[3.] On a la m\^eme formule par sym\'etrie si $q_2=1$ et $q_1=p_j>1$.
\item[4.] Enfin, si $q_1=q_2=1$, alors
$$
\begin{aligned}
(d_\ell c)^{(N+1)}_{1,1,p_1\dots\check{j}\dots p_n}&(Y_1\pt Y_2\pt X_1\pt\dots\hat\jmath\dots\pt X_n)=\\
=&\varepsilon\left(\begin{smallmatrix}y_1y_2x_1\dots\hat\jmath\dots x_n\\ y_1x_1\dots y_2\dots x_n\end{smallmatrix}\right)
Y_1\ppt c^{(N)}_{p_1\dots p_N}(X_1\pt\dots\pt Y_2\pt\dots\pt X_n)\\
&+\varepsilon\left(\begin{smallmatrix}y_1y_2x_1\dots\hat\jmath\dots x_n\\ y_2x_1\dots y_1\dots x_n\end{smallmatrix}\right)
Y_2\ppt c^{(N)}_{p_1\dots p_N}(X_1\pt\dots\pt Y_1\pt\dots\pt X_n)\\
&-(-1)^{\sum_{i<j}\deg(X_i)}\varepsilon\left(\begin{smallmatrix}y_1y_2x_1\dots\hat\jmath\dots x_n\\ x_1\dots y_1y_2\dots x_n\end{smallmatrix}\right) c^{(N)}_{p_1\dots p_n}(X_1\pt\dots\pt\ell(Y_1\pt Y_2)\pt\dots\pt X_n).
\end{aligned}
$$
\end{itemize}

On a une cohomologie : $(d_\ell+d_m)\circ(d_\ell+d_m)=0$.\\
\end{prop}


\

\end{document}